\documentclass[11pt]{article}
\usepackage{amssymb,amsmath,amsthm,mathrsfs}
\usepackage{times,a4wide}
\usepackage{lineno}
\usepackage{cite}
\smallskip

\newtheorem{theorem}{Theorem}[section]

\newtheorem{proposition}[theorem]{Proposition}
\newtheorem{lemma}[theorem]{Lemma}
\newtheorem{definition}[theorem]{Definition}

\newtheorem{remark}[theorem]{Remark}
\numberwithin{equation}{section}

\title{{\bf Multiple solutions for a fractional Choquard problem with slightly subcritical exponents on bounded domains}}

\author{
	{\sc Marco Ghimenti$^{1}$  \thanks{Email: marco.ghimenti@unipi.it}}, 	\,
	{\sc  Min Liu$^{2}$ \thanks{Email: minliu@lnnu.edu.cn; Supported by NSFC (12101282), LJKZ0967 and 2021BSL004; Corresponding author}}   \,	and	
	{\sc Zhongwei Tang$^{3}$  \thanks{Email: tangzw@bnu.edu.cn; Supported by NSFC (12071036,11671331)}}\\		
	{\small $^{1}$ Dipartimento di Matematica, Universit\`a di Pisa, Largo B. Pontecorvo 5, 56126 Pisa, Italy}\\
	{\small $^{2}$ School of Mathematics, Liaoning Normal University, Dalian, 116029, China}\\
	{\small $^{3}$ School of Mathematical Sciences, Beijing Normal University, Beijing, 100875, China}
}

\date{}
\makeatletter

\begin{document}

\maketitle

\begin{abstract}
This paper is devoted to study a fractional Choquard problem with slightly subcritical exponents on bounded domains. When the exponent of the convolution type nonlinearity tends to the fractional critical one in the sense of Hardy-Littlewood-Sobolev inequality, we obtain the existence of multiple positive solutions via Lusternik-Schnirelmann category and nonlocal global compactness. Moreover, we prove that the topology of the domain furnishes a lower bound for the number of positive solutions.
\end{abstract}

{\bf Keywords}: Fractional Choquard problem; Slightly subcritical exponent; Category; Nonlocal global compactness.

{\bf AMS Subject Classification}: 35J66; 35Q40; 35R11.

\section{Introduction and a main result}

In the present paper, we are concerned with the following fractional Choquard problem with slightly subcritical exponents on bounded domains:
\begin{equation}\label{target}
    \begin{cases}
    (-\Delta)^s u+\lambda u=\left(\displaystyle\int_\Omega\frac{u^{p_\varepsilon}(y)}{|x-y|^\sigma}\text{d}y\right)u^{p_\varepsilon-1} \quad &\text{in }\Omega\\
    u=0 \quad &\text{in }\mathbb R^N\setminus\Omega\\
    u>0 \quad &\text{in }\Omega
    \end{cases},
\end{equation}
where $s\in(0,1),\ \lambda\ge0,\ N>4s,\ \sigma\in(0,N),\ \varepsilon>0,\ p_\varepsilon:=2^*_{\sigma,s}-\varepsilon,\ 2^*_{\sigma,s}:=\frac{2N-\sigma}{N-2s}$ is the fractional critical exponent in the sense of Hardy-Littlewood-Sobolev inequality, $\Omega$ is a bounded domain in $\mathbb R^N$ with Lipschitz boundary, and $(-\Delta)^s$ is the fractional Laplacian defined by
\begin{equation}\label{fld}
(-\Delta)^su(x)=C_{N,s}\lim_{\rho\to0}\int_{\mathbb{R}^{N}\setminus B_\rho(x)} \frac{u(x)-u(y)}{|x-y|^{N+2s}}\text{d}y
\end{equation}
with $C_{N,s}:=\big(\int_{\mathbb{R}^N}\frac{1-\cos\zeta_1}{|\zeta|^{N+2s}}\text{d}\zeta\big)^{-1}$ being a dimensional constant and $B_\rho(x)$ being an open ball centered at $x$ with radius $\rho$.
As $\varepsilon\to0,$ namely, $p_\varepsilon\to 2^*_{\sigma,s}$, we prove that the problem \eqref{target} possesses at least $cat_{\bar\Omega}(\bar\Omega)+1$ solutions if $\Omega$ is not contractible. Here $cat_{\bar\Omega}(\bar\Omega)$ denotes the Lusternik-Schnirelmann category of $\Omega$.

On this type of problem in history, we can go back to Bahri and Coron \cite{Bahri-Coron}, Benci and Cerami \cite{Benci-Cerami} and Benci, Cerami and Passaseo \cite{Benci-Cerami-Passaseo}. In the celebrated work \cite{Benci-Cerami-Passaseo}, the authors considered
\begin{equation}\label{Benci-Cerami-Passaseo}
\begin{cases}
-\Delta u+\lambda u=u^{p-1} \quad &\text{in }\Omega\\
u=0 \quad &\text{on } \partial\Omega\\
u>0 \quad &\text{in }\Omega
\end{cases},
\end{equation}
where $N\ge 3$, $\lambda\ge0$, $p\in(2,2^*:=\frac{2N}{N-2})$ and $\Omega\subset\mathbb R^N$ is a smooth bounded domain. They showed to us how the number of solutions for \eqref{Benci-Cerami-Passaseo} is affected by the topology of $\Omega$, and the nonlinearity acts strongly on the problem \eqref{Benci-Cerami-Passaseo} if the domain $\Omega$ is perturbed by cutting off or adding pieces, small in some sense. To be precise, they applied variational methods to prove that there exists a function $\bar p:[0,+\infty)\to(2,2^*)$ such that for every $p\in[\bar p(\lambda),2^*)$, the problem \eqref{Benci-Cerami-Passaseo} has at least $cat_{\bar\Omega}(\bar\Omega)$ distinct solutions. Moreover, they also showed that the number of solutions is greater than $cat_{\bar\Omega}(\bar\Omega)+1$ when the domain is not contractible. More than a decade later, Benci, Bonanno and Micheletti \cite{Benci-Bonanno-Micheletti} extended this kind of result to a nonlinear elliptic problem on a Riemannian manifold and proved that the number of solutions depends on the topological properties of the manifold. In the same spirit, Siciliano \cite{Siciliano} investigated the existence of multiple positive solutions to a Schr\"odinger-Poisson-Slater system 
\begin{equation}\label{Siciliano}
\begin{cases}
-\Delta u+u+\lambda \phi u=|u|^{p-2}u \quad &\text{in }\Omega\\
-\Delta \phi=u^2 \quad &\text{in }\Omega\\
u=\phi=0 \quad &\text{on } \partial\Omega
\end{cases},
\end{equation}
and showed that the number of positive solutions to \eqref{Siciliano} is bounded from below by $cat_{\bar\Omega}(\bar\Omega)+1$ under appropriate assumptions. Afterwards, this thought was also used to study nonlinear complex equations with magnetic field. For instance, Alves, Figueiredo and Furtado \cite{Alves} established an existence result of multiple solutions for a nonlinear complex Schr\"odinger equation with magnetic field when the parameter $\lambda$ has large values. Alves, Figueiredo and Yang \cite{Alves-2} also obtained multiple solutions of a nonlinear complex Choquard equation with magnetic field and related the number of solutions with the topology of the set where the potential attains its minimum value. Subsequently, Xiang, R$\breve a$dulescu and Zhang \cite{Xiang} considered a fractional Choquard–Kirchhoff type problem involving an external magnetic potential and got the existence of nontrivial radial solutions in non-degenerate and degenerate cases, respectively. Ji and R$\breve a$dulescu \cite{Ji} studied the multiplicity of multi-bump solutions for a nonlinear magnetic Choquard equation with deepening potential well when the zero set of the potential has several isolated connected components. 

We also want to mention some recent work about multiplicity  for Choquard equations dealing with concentration properties. Yang and Zhao \cite{Yang} applied penalization techniques and Ljusternik-Schnirelmann theory to investigate the relation between the number of positive solutions and the topology of the set where the potential attains its minimum values for a singularly perturbed fractional Choquard equation, and showed that these positive solutions locate at the minimum point of the potential. In addition, Cingolani and Tanaka \cite{Cingolani} studied the multiplicity and concentration of positive single-peak solutions for a nonlinear Choquard equation, and they used relative cup-length to estimate the topological changes and then related the number of positive solutions to the topology of the critical set of the potential.

Recently, based on the work of Gao and Yang \cite{Gao-Yang}, Ghimenti and Pagliardini \cite{Ghimenti-Pagliardini} used the concentration properties of the Talenti-Aubin functions to balance the nonlocal effect of the nonlinerity and studied a Choquard equation in a bounded domain
\begin{equation}\label{Ghimenti-Pagliardini}
    \begin{cases}
    -\Delta u+\lambda u=\left(\displaystyle\int_\Omega\frac{u^{q_\varepsilon}(y)}{|x-y|^\sigma}\text{d}y\right)u^{q_\varepsilon-1} \quad &\text{in }\Omega\\
    u=0 \quad &\text{on }\partial\Omega\\
    u>0 \quad &\text{in }\Omega
\end{cases},
\end{equation}
where $\Omega$ is a regular bounded domain in $\mathbb R^N$, $N>3$, $\sigma\in(0,N)$, $\lambda\ge0$, $\varepsilon>0$ and $q_\varepsilon=2^*_\sigma-\varepsilon$ with $2^*_\sigma:=\frac{2N-\sigma}{N-2}$. They also achieved a similar result, which is the existence of at least $cat_{\bar\Omega}(\bar\Omega)+1$ solutions for \eqref{Ghimenti-Pagliardini} when the exponent $q_\varepsilon$ is close to $2^*_\sigma$ and the domain is not contractible.

The aim of this paper is to investigate the fractional counterpart of the above problem \eqref{Ghimenti-Pagliardini}. Indeed, according to the large amount of literatures on fractional Laplacian appearing in the last years, it is very natural to ask whether a similar result holds for the corresponding problem to \eqref{Ghimenti-Pagliardini} with the fractional Laplacian. Therefore, we would like to study the slightly critical problem \eqref{target} with two nonlocal framework caused by the fractional Laplacian and the convolution term, and finally implement a multiplicity result. The nonlocality of fractional Laplacian makes the discussion different from that in \cite{Ghimenti-Pagliardini}, especially the proof of the crucial nonlocal splitting lemma.

It is known that if $\Omega$ is a bounded starshaped domain and the nonlinearity is of critical growth, namely,  $p_\varepsilon$ is replaced by $2^*_{\sigma,s}$, then the problem \eqref{target} has no solution according to the Pohozaev identity (refer to \cite{Ros-Oton})
$$-s\lambda\int_\Omega u^2\text{d}x=\frac{\Gamma(s+1)^2}{2}\int_{\partial\Omega}\left(\frac{u}{\delta^s}\right)^2(x\cdot\nu)\text{d}S,\quad \text{where } \delta=\text{dist}(x,\partial\Omega).$$
However, if the critical exponent is perturbed slightly, namely, $p_\varepsilon=2^*_{\sigma,s}-\varepsilon$, then the existence of solution can be shown and even multiple solutions exist. That is what we study in this paper. We also want to mention \cite{Mukherjee-Sreenadh} in which equation \eqref{target} with critical nonlinearity (that is $\varepsilon=0$) is studied on a bounded domain and the analogous of Br\'ezis-Nirenberg problem is proved. The main difference between the two problems is that we use the perturbation of the exponent to obtain the multiplicity result, while in the paper of Mukherjee and Sreenadh is the parameter $\lambda$ in front of the linear term which plays a key role.

In order to complete the study for \eqref{target}, we shall apply variational methods, nonlocal global compactness, Lusternik-Schnirelmann category theory and the technique introduced by \cite{Benci-Cerami,Benci-Cerami-Passaseo} to make a comparison between the category of some sublevel set of the functional and the category of domain $\Omega$. Hence we consider two limit problems and make a careful study on the behavior of the associated functionals and their related minimal levels, where the concentration property of the Talenti-Aubin functions is applied for more than once. About the nonlocal global compactness, we would like to mention a recent work by He and R$\breve a$dulescu \cite{He}, in which they obtained a very useful property (see Lemma \ref{He}) for the (PS)-sequence to the limit problem \eqref{limit1}. By means of this property, we establish a crucial splitting lemma (Lemma \ref{split}) to the second limit problem \eqref{limit2}. Together with the barycenter map \eqref{barycenter}, we can find the sublevels of the functional with barycenter not away from $\Omega$ and category greater than $cat_{\bar\Omega}(\bar\Omega)$, and then get the existence of at least $cat_{\bar\Omega}(\bar\Omega)$ solutions as $\varepsilon\to0$ according to the Lusternik-Schnirelmann category theory. Since the fractional Laplacian makes nonlocal effect, the barycenter map here is different from that in Ghimenti and Pagliardini \cite{Ghimenti-Pagliardini}. In order to give the existence of another solution when $\Omega$ is not contractible, we construct a compact and contractible set $T_\varepsilon$ containing only positive functions and being a subset of the Nehari manifold. 

Apart from the case of $p\to2^*$, the works in \cite{Benci-Cerami,Benci-Cerami-Passaseo} also treat the situation of $\lambda\to+\infty$. Although in this paper we don't explore the second situation in detail, our result still tells that the domain topology gives a lower bound on the number of solutions to \eqref{target} when the parameter $\lambda$ is very large.

Now we state our main result, which gives an affirmative answer on the possibility of extending the result in Ghimenti and Pagliardini \cite{Ghimenti-Pagliardini} to the fractional case.

\begin{theorem}\label{t}
	Assume that $s\in(0,1)$, $\lambda\ge0$, $N>4s$ and $\sigma\in(0,N)$. Then there exists $\hat\varepsilon>0$ such that for each $\varepsilon\in(0,\hat\varepsilon]$, the problem \eqref{target} has at least $cat_{\bar\Omega}(\bar\Omega)$ low energy solutions. Moreover, if $\Omega$ is not contractible, there is another solution with higher energy.
\end{theorem}

\begin{remark}
	The assumption $N>4s$ is needed to estimate some integrals related to the bubbles of the limit problem, such as \eqref{1}, \eqref{2} and \eqref{3}. Indeed, it is \eqref{3} that requires $N>4s$, while for \eqref{1} and \eqref{2}, $N>2s$ is enough. 
\end{remark}

\begin{remark}
	If $cat_{\bar\Omega}(\bar\Omega)=1$, then it's easy to achieve the existence of one solution. Indeed, it can be obtained in a simplier way by the famous Mountain Pass Theorem. Additionally, if $\Omega$ is not contractible, we can obtain $cat_{\bar\Omega}(\bar\Omega)$ low energy solutions and another high energy solution. Therefore, we assume $cat_{\bar\Omega}(\bar\Omega)>1$ in what follows.
\end{remark}

Our paper in the following is organized like this: Section $2$ contains the working spaces, some essential definitions, several classical results and useful relations. Section $3$ gives the variational structure of this problem, the functional setting and some preliminary results. Section $4$ serves two limit problems and especially a nonlocal splitting lemma, which play a key role in the arguments of the main theorem. Section $5$ contributes to the proof of our main result.

\section{Notations and preliminaries}

For $s\in(0,1)$, we denote $2^*_s:=\frac{2N}{N-2s}$ as the fractional critical Sobolev exponent, and define two Hilbert spaces:
\begin{equation*}
    H^s(\mathbb{R}^N)=\left\{u\in L^2(\mathbb{R}^N):\frac{|u(x)-u(y)|}{|x-y|^{\frac{N+2s}{2}}}\in L^2(\mathbb{R}^N\times\mathbb{R}^N)\right\}
\end{equation*}
with the inner product
\begin{equation*}
    (u,v)_{H^s(\mathbb{R}^N)}=\int_{\mathbb{R}^{2N}}\frac{(u(x)-u(y))(v(x)-v(y))}{|x-y|^{N+2s}}\text{d}x\text{d}y+\int_{\mathbb{R}^N}uv\text{d}x
\end{equation*}
and the norm
\begin{equation*}
   \|u\|_{H^s(\mathbb{R}^N)}=\left(\int_{\mathbb{R}^{2N}}\frac{|u(x)-u(y)|^2}{|x-y|^{N+2s}}\text{d}x\text{d}y+\int_{\mathbb{R}^N}u^2\text{d}x\right)^\frac{1}{2},
\end{equation*}
and
\begin{equation*}
   D^{s,2}(\mathbb{R}^N)=\left\{u\in L^{2^*_s}(\mathbb{R}^N):\frac{|u(x)-u(y)|}{|x-y|^{\frac{N+2s}{2}}}\in L^2(\mathbb{R}^N\times\mathbb{R}^N)\right\}
\end{equation*}
with the inner product
\begin{equation*}
   (u,v)_{D^{s,2}(\mathbb{R}^N)}=\int_{\mathbb{R}^{2N}}\frac{(u(x)-u(y))(v(x)-v(y))}{|x-y|^{N+2s}}\text{d}x\text{d}y
\end{equation*}
and the norm
\begin{equation*}
   \|u\|_{D^{s,2}(\mathbb{R}^N)}=\left(\int_{\mathbb{R}^{2N}}\frac{|u(x)-u(y)|^2}{|x-y|^{N+2s}}\text{d}x\text{d}y\right)^\frac{1}{2}.
\end{equation*}
Here we have used the fractional Sobolev inequality (see Theorem 1.1 in Cotsiolis et al \cite{Cotsiolis} and Proposition 3.6 in Di Nezza et al \cite{Nezza})
\begin{equation}\label{FSI}
	\|u\|^2_{L^{2^*_s}(\mathbb{R}^N)}\le	S_{N,s}\int_{\mathbb R^{2N}}\frac{|u(x)-u(y)|^2}{|x-y|^{N+2s}}\text{d}x\text{d}y, \quad\forall u\in C^\infty_0(\mathbb R^N),
\end{equation}
where $S_{N,s}:=\frac{1}{2}(4\pi)^{-s}\frac{\Gamma(\frac{N-2s}{2})}{\Gamma(\frac{N+2s}{2})}\Big[\frac{\Gamma(N)}{\Gamma(\frac{N}{2})}\Big]^{\frac{2s}{N}}C_{N,s}$ is achieved if and only if $$u(x)=A(\alpha^2+(x-x_0)^2)^{-\frac{N-2s}{2}},\quad x\in\mathbb{R}^N$$ for some fixed constant $A\in \mathbb{R}$ and parameters $\alpha\in\mathbb R\setminus\{0\}, \ x_0\in\mathbb{R}^N$. Meanwhile, \eqref{FSI} means that $D^{s,2}(\mathbb{R}^N)$ is continuously embedded in $L^{2^*_s}(\mathbb R^N)$.

Let 
\begin{equation*}
    H^s_0(\Omega)=\left\{u\in H^s(\mathbb{R}^N):u=0\text{ a.e. in }\mathbb R^N\setminus\Omega\right\}
\end{equation*}
and
\begin{equation*}
    D^{s,2}_0(\Omega)=\left\{u\in D^{s,2}(\mathbb{R}^N):u=0\text{ a.e. in }\mathbb R^N\setminus\Omega\right\}.
\end{equation*}
Since $\Omega$ is a bounded domain with Lipschitz boundary, it follows from Corollary 5.5 in Di Nezza \cite{Nezza} that $H^s_0(\Omega)$ can be regarded as the closure of $C^\infty_0(\Omega)$ in $H^s(\mathbb{R}^N)$ and $D^{s,2}_0(\Omega)$ is the closure of $C^\infty_0(\Omega)$ in $D^{s,2}(\mathbb{R}^N)$. Moreover, $H^s_0(\Omega)=D^{s,2}_0(\Omega).$ According to Theorems 6.5 and 7.2 in Di Nezza et al \cite{Nezza}, $H^s_0(\Omega)$ is continuously embedded in $L^q(\Omega)$ for $q\in[1,2^*_s]$, and is compactly embedded in $L^q(\Omega)$ for $q\in[1,2^*_s)$. 
Notice that $\lambda\ge0$, we can choose  
\begin{equation*}
\|u\|_\lambda:=\left(\int_{\Upsilon}\frac{|u(x)-u(y)|^2}{|x-y|^{N+2s}}\text{d}x\text{d}y+\lambda\int_{\Omega}u^2\text{d}x\right)^\frac{1}{2}
\end{equation*}
as an equivalent norm in $H^s_0(\Omega)$, where $\Upsilon:=\mathbb R^{2N}\setminus (\Omega^c\times\Omega^c)$ with $\Omega^c:=\mathbb R^N\setminus \Omega$.



Recall the well-known Lion's Lemma, which will be useful for showing the achievement of the least energy on Nehari manifold.

\begin{lemma}\label{Lions}
	(Felmer et al \cite{Felmer})
	Let $s\in (0,1)$ and $r>0$. If $\{u_n\}_{n\in\mathbb N}$ is bounded in $H^s(\mathbb{R}^N)$ and 
	\begin{equation*}
	\sup_{y\in\mathbb{R}^N}\int_{B_r(y)}|u_n(x)|^2\text{d}x\to0\quad\text{as }\ n\to\infty,
	\end{equation*}
	then $u_n\to 0$ in $L^q(\mathbb{R}^N)$ for $q\in(2,2_s^*)$.
\end{lemma}

Now we recall some information on Lusternik-Schnirelmann category.

\begin{definition}
	Let $A$ be a closed subset of a topological space $X$. The category of $A$ in $X$, denoted by $cat_X(A)$, is the least integer $k$ such that $A\subseteq A_1\cup A_2\cup\cdots\cup A_k$ with $A_i$ closed and contractible in $X$ for each $i=1,2,\cdots,k.$
\end{definition}

We set $cat_X(\varnothing)=0$, and $cat_X(A)=+\infty$ if there is no integer with the above property. Without confusion, we write $cat(X)$ for $cat_X(X)$.

\begin{remark}\label{BC}
	(Benci et al \cite{Benci-Cerami})
	Let $X$ and $Y$ be two topological spaces. If $f:X\to Y$ and $g:Y\to X$ are continuous operators such that $g\circ f$ is homotopic to the identity on $X$, then $cat(X)\le cat(Y)$.
\end{remark}

\begin{proposition}\label{BCP}
	(Benci et al \cite{Benci-Cerami-Passaseo})
	Let $J$ be a $C^1$ real functional on a complete $C^{1,1}$ Banach manifold $M$. If $J$ is bounded from below and satisfies Palais-Smale condition on $M$, then $J$ has at least $cat(J^d)$ critical points in $J^d$, where $J^d:=\{u\in M:J(u)\le d\}$. Moreover, if $M$ is contractible and $cat(J^d)>1$, then there is at least one critical point $u\notin J^d$.
\end{proposition}

We also remind the famous Hardy-Littlewood-Sobolev inequality and some results linked to it.

\begin{lemma}\label{HLS}
	(Lieb et al \cite{Lieb-Loss})
	Let $r,t>1$ and $0<\sigma<N$ with $\frac{1}{r}+\frac{\sigma}{N}+\frac{1}{t}=2$, $f\in L^r(\mathbb R^N)$ and $g\in L^t(\mathbb R^N)$. Then there exists a constant $C_{N,\sigma,r}$ independent of $f$ and $g$ such that 
	\begin{equation*}
	\bigg|\int_{\mathbb R^{2N}}\frac{f(x)g(y)}{|x-y|^\sigma}\text{d}x\text{d}y\bigg|\le C_{N,\sigma,r}\|f\|_{L^r(\mathbb R^N)}\|g\|_{L^t(\mathbb R^N)}.
	\end{equation*}
	Moreover, this equality holds if and only if $r=t=\frac{2N}{2N-\sigma}$, $g=Cf$ ($C$ is a constant) and 
	\begin{equation*}
	f(x)=B(\alpha^2+|x-x_0|^2)^{-\frac{2N-\sigma}{2}},\quad x\in\mathbb R^N
	\end{equation*}
	for some $B\in\mathbb C$, $\alpha\in\mathbb R\setminus\{0\}$ and $x_0\in\mathbb R^N$.
\end{lemma}

Due to Lemma \ref{HLS}, embedding result and the equivalence of the norms $\|\cdot\|$ and $\|\cdot\|_\lambda$, we have the following result.

\begin{remark}\label{riesz<norm}
	For any $u\in H^s_0(\Omega)$, there exists a constant $C>0$ such that 
	\begin{equation*}
	\int_{\Omega\times\Omega}\frac{|u(x)|^{p_\varepsilon}|u(y)|^{p_\varepsilon}}{|x-y|^\sigma}\text{d}x\text{d}y\le C_{N,\sigma,s,\varepsilon}\Big\|u\Big\|^{2p_\varepsilon}_{L^{\frac{2Np_\varepsilon}{2N-\sigma}}(\Omega)}\le C\|u\|^{2p_\varepsilon}_\lambda.
	\end{equation*}
\end{remark}

Using again the Hardy-Littlewood-Sobolev inequality and embedding result, we also get that for all $u\in D^{s,2}(\mathbb R^N)$,
\begin{equation}\label{critical<}
    \bigg(\int_{\mathbb R^{2N}}\frac{|u(x)|^{2^*_{\sigma,s}}|u(y)|^{2^*_{\sigma,s}}}{|x-y|^\sigma}\text{d}x\text{d}y\bigg)^{\frac{1}{2^*_{\sigma,s}}}\le C_{N,\sigma,s}\|u\|^2_{L^{2^*_s}(\mathbb R^N)}.
\end{equation}
Define 
\begin{equation}\label{S1}
    S^s_{H,L}:=\inf_{u\in D^{s,2}(\mathbb R^N)\setminus\{0\}}\frac{\displaystyle\int_{\mathbb R^{2N}}\frac{|u(x)-u(y)|^2}{|x-y|^{N+2s}}\text{d}x\text{d}y}{\bigg(\displaystyle\int_{\mathbb R^{2N}}\frac{|u(x)|^{2^*_{\sigma,s}}|u(y)|^{2^*_{\sigma,s}}}{|x-y|^\sigma}\text{d}x\text{d}y\bigg)^{\frac{1}{2^*_{\sigma,s}}}}
\end{equation}
and it follows from Lemma \ref{HLS}, \eqref{FSI} and \eqref{critical<} that $S^s_{H,L}$ is achieved if and only if 
\begin{equation}\label{U}
    u(x)=U_{R,a}(x):=A\left(\frac{R}{1+R^2|x-a|^2}\right)^{\frac{N-2s}{2}},\quad x\in\mathbb R^N,
\end{equation}
where $A>0$ is some fixed constant, $R>0$ and $a\in\mathbb R^N$ are parameters. Absolutely, $S^s_{H,L}=(S_{N,s}C_{N,\sigma,s})^{-1}$.
We can also define 
\begin{equation}\label{SHL}
	S^{s,\Omega}_{H,L}:=\inf_{u\in D^{s,2}_0(\Omega)\setminus\{0\}}\frac{\displaystyle\int_{\Upsilon}\frac{|u(x)-u(y)|^2}{|x-y|^{N+2s}}\text{d}x\text{d}y}{\bigg(\displaystyle\int_{\Omega\times\Omega}\frac{|u(x)|^{2^*_{\sigma,s}}|u(y)|^{2^*_{\sigma,s}}}{|x-y|^\sigma}\text{d}x\text{d}y\bigg)^{\frac{1}{2^*_{\sigma,s}}}}.
\end{equation}
By means of Lemma 2.2 in Mukherjee et al \cite{ Mukherjee-Sreenadh}, $S^{s,\Omega}_{H,L}=S^s_{H,L}$ and $S^{s,\Omega}_{H,L}$ is never achieved unless $\Omega=\mathbb R^N$.

In the following arguments, we denote
\begin{equation*}
    \|u\|_{s,2}:=\|u\|_{D^{s,2}(\mathbb{R}^N)},\quad |u|_q:=\|u\|_{L^q(\mathbb{R}^N)},
\end{equation*}
\begin{equation*}
    u^+:=\max\{0,u\}, \quad u^-:=\min\{0,u\},\quad \mathbb R^+:=(0,+\infty).
\end{equation*}	
For simplicity and without destruction, we drop the constant $C_{N,s}$ in the definition of $(-\Delta)^s$, and we shall use $C$ to represent various positive constants, which may be different in different places.

\section{Variational setting}

To study the problem \eqref{target}, we consider the associated functional $J_\varepsilon:H^s_0(\Omega)\to\mathbb R$ given by
\begin{equation*}
    J_\varepsilon(u):=\frac{1}{2}\|u\|^2_\lambda-\frac{1}{2p_\varepsilon}\int_{\Omega\times\Omega}\frac{(u^+(x))^{p_\varepsilon}(u^+(y))^{p_\varepsilon}}{|x-y|^\sigma}\text{d}x\text{d}y,
\end{equation*}
and the associated Nehari manifold 
\begin{equation*}
    \mathcal N_\varepsilon:=\left\{u\in H^s_0(\Omega)\setminus\{0\}:\left<J'_\varepsilon(u),u\right>=0\right\}.
\end{equation*}
We set the least energy by
\begin{equation*}
    \vartheta_\varepsilon:=\inf_{u\in\mathcal N_\varepsilon}J_\varepsilon(u).
\end{equation*}
It is easy to verify that $J_\varepsilon\in C^1(H^s_0(\Omega),\mathbb R)$, since for any $v\in H^s_0(\Omega)$,
\begin{align*}
    \left<J'_\varepsilon(u),v\right>=&\int_\Upsilon\frac{(u(x)-u(y))(v(x)-v(y))}{|x-y|^{N+2s}}\text{d}x\text{d}y+\lambda\int_\Omega uv\text{d}x\\
    &-\int_{\Omega\times\Omega}\frac{(u^+(y))^{p_\varepsilon}(u^+(x))^{p_\varepsilon-1}v(x)}{|x-y|^\sigma}\text{d}x\text{d}y.
\end{align*}

Now we establish some preliminary results.

\begin{proposition}\label{equivalence}
	The nonzero critical points of $J_\varepsilon$ coincide with the solutions of the problem \eqref{target}.
\end{proposition}

\begin{proof}[\bf Proof]
    It is clear that the nonzero critical points of $J_\varepsilon$ are solutions of
    \begin{equation}\label{new}
    \begin{cases}
    (-\Delta)^s u+\lambda u=\left(\displaystyle\int_\Omega\frac{(u^+(y))^{p_\varepsilon}}{|x-y|^\sigma}\text{d}y\right)(u^+)^{p_\varepsilon-1}\quad &\text{in }\Omega\\
    u=0 \quad &\text{in }\mathbb R^N\setminus\Omega\\
    u\neq 0\quad &\text{in }\Omega
    \end{cases}.
    \end{equation}
    We claim that \eqref{new} is equivalent to \eqref{target}. 
    
    On one hand, if $u\in H^s_0(\Omega)$ is a solution of \eqref{target}, then $u=u^+$, and thus $u$ is a solution of \eqref{new}.
    
    On the other hand, if $u\in H^s_0(\Omega)$ is a solution of \eqref{new}, then $u$ is a critical point of $J_\varepsilon$. It follows from 
    \begin{align*}
    \int_{\Upsilon}\frac{|u^-(x)-u^-(y)|^2}{|x-y|^{N+2s}}\text{d}x\text{d}y=&\int_{\{u>0\}\times\{u<0\}}\frac{|-u(y)|^2}{|x-y|^{N+2s}}\text{d}x\text{d}y\\
    &+\int_{\{u<0\}\times\{u>0\}}\frac{|u(x)|^2}{|x-y|^{N+2s}}\text{d}x\text{d}y\\
    &+\int_{\{u<0\}\times\{u<0\}}\frac{|u(x)-u(y)|^2}{|x-y|^{N+2s}}\text{d}x\text{d}y\\
    \le&\int_{\{u>0\}\times\{u<0\}}\frac{|u(x)-u(y)|^2}{|x-y|^{N+2s}}\text{d}x\text{d}y\\
    &+\int_{\{u<0\}\times\{u>0\}}\frac{|u(x)-u(y)|^2}{|x-y|^{N+2s}}\text{d}x\text{d}y\\
    &+\int_{\{u<0\}\times\{u<0\}}\frac{|u(x)-u(y)|^2}{|x-y|^{N+2s}}\text{d}x\text{d}y\\
    \le&\int_{\Upsilon}\frac{|u(x)-u(y)|^2}{|x-y|^{N+2s}}\text{d}x\text{d}y
    \end{align*}
    that $u^-\in H^s_0(\Omega)$. Now we use $u^-$ as a test function and obtain that
    \begin{align*}
    0=\left<J'_\varepsilon(u),u^-\right>=&\int_\Upsilon\frac{(u(x)-u(y))(u^-(x)-u^-(y))}{|x-y|^{N+2s}}\text{d}x\text{d}y+\lambda\int_\Omega(u^-)^2\text{d}x\\
    &+\int_{\Omega\times\Omega}\frac{(u^+(y))^{p_\varepsilon}(u^+(x))^{p_\varepsilon-1}u^-(x)}{|x-y|^\sigma}\text{d}x\text{d}y\\
    =&\int_{\{u<0\}\times\{u<0\}}\frac{(u^-(x)-u^-(y))^2}{|x-y|^{N+2s}}\text{d}x\text{d}y\\
    &+\int_{\{u>0\}\times\{u<0\}}\frac{(u^+(x)-u^-(y))(-u^-(y))}{|x-y|^{N+2s}}\text{d}x\text{d}y\\
    &+\int_{\{u<0\}\times\{u>0\}}\frac{(u^-(x)-u^+(y))u^-(x)}{|x-y|^{N+2s}}\text{d}x\text{d}y+\lambda\int_\Omega(u^-)^2\text{d}x\\
    \ge&\int_\Upsilon\frac{(u^-(x)-u^-(y))^2}{|x-y|^{N+2s}}\text{d}x\text{d}y+\lambda\int_\Omega(u^-)^2\text{d}x\\
    =&\|u^-\|^2_\lambda,
    \end{align*}
    which implies $u^-=0$. Hence $u=u^+\ge0$. By Theorem 3.2 in d'Avenia et al \cite{Avenia}, we get $u\in C^{0,\mu}(\mathbb{R}^N)$ with $\mu\in(0,1)$. Thus $u\in L^\infty(\mathbb{R}^N)$. Due to Theorem 1 and its proof in Du Plessis \cite{Plessis}, we have $|x|^{-\sigma}*u^{p_\varepsilon}\in C^{0,\mu+\sigma}(\mathbb{R}^N)$ if $\mu+\sigma\in(0,1)$, and $|x|^{-\sigma}*u^{p_\varepsilon}\in C^{0,1}(\mathbb{R}^N)$ if $\mu+\sigma\in[1,N+1)$. As Theorem 1.4 in Felmer et al \cite{Felmer}, we obtain that $u\in C^{0,2s+\mu}(\mathbb{R}^N)$ if $2s+\mu\le1$, and $u\in C^{1,2s+\mu-1}(\mathbb{R}^N)$ if $2s+\mu>1$. This regularity makes sure \eqref{fld} hold for $u$ in the pointwise sense. Suppose $u(x_0)=0$ for some $x_0\in\Omega$, then the equation in \eqref{new} deduces $\left((-\Delta)^su\right)(x_0)=0$. By using \eqref{fld} again, we get $u\equiv0$, which is impossible. Hence $u>0$ in $\Omega$. Concequently, $u$ is also a solution of \eqref{target}.  
\end{proof}

By Remark \ref{riesz<norm}, we can get the common property of Nehari manifold, that is, there exists a constant $C>0$ such that $\|u\|_\lambda>C$ and $J_\varepsilon(u)>C$ for all $u\in\mathcal N_\varepsilon$. In addition, we also see that for any $0\neq u\in H^s_0(\Omega)$, there is a unique constant $t_\varepsilon(u)>0$ such that $t_\varepsilon(u)u\in\mathcal N_\varepsilon$. Based on these information and a standard proof, we have the following result.

\begin{lemma}\label{energy}
	There holds
	$$0<\vartheta_\varepsilon=\inf_{u\in H^s_0(\Omega)\setminus\{0\}}\max_{t\ge0}J_\varepsilon(tu)=\inf_{\gamma\in\Gamma_\varepsilon}\max_{t\in[0,1]}J_\varepsilon(\gamma(t)),$$
	where
	\begin{equation*}
	\Gamma_\varepsilon:=\left\{\gamma\in C\left([0,1],H^s_0(\Omega)\right):\gamma(0)=0,\ J_\varepsilon(\gamma(1))<0\right\}.
	\end{equation*}
\end{lemma}

From Lemma \ref{Lions}, Lemma \ref{energy}, and a similar proof to that of Theorem 1 in Moroz et al \cite{Moroz-JFA} or Lemma 3.5 in Liu et al \cite{Liu-DCDS}, the following result is true.

\begin{lemma}\label{groundstate}
	$\vartheta_\varepsilon$ is achieved by a function $u_\varepsilon\in\mathcal N_\varepsilon$, that is $\vartheta_\varepsilon=J_\varepsilon(u_\varepsilon)$.
\end{lemma}

Finally, we review Palais-Smale sequence simplified by (PS)-sequence.

\begin{definition}
	A sequence $\{u_n\}_{n\in\mathbb N}$ is called a  (PS)-sequence of $J_\varepsilon$, if 
	$$\{J_\varepsilon(u_n)\}_{n\in\mathbb N} \text{ is bounded in } \mathbb R \text{ and } J'_\varepsilon(u_n)\to0 \text{ in } H^{-s}_0(\Omega),$$
	where $H^{-s}_0(\Omega)$ signifies the dual space of $H^s_0(\Omega)$. If every (PS)-sequence of $J_\varepsilon$ has a convergent subsequence, then we say $J_\varepsilon$ satisfies the (PS)-condition on $H^s_0(\Omega)$. 
\end{definition}

In fact, $J_\varepsilon$ does satisfy the (PS)-condition on $H^s_0(\Omega)$ globally. 

\begin{lemma}\label{PS}
	If $\{u_n\}_{n\in\mathbb N}\subset\mathcal N_\varepsilon$ is a (PS)-sequence of the constrained functional $J_\varepsilon|_{\mathcal N_\varepsilon}(u)=\frac{p_\varepsilon-1}{2p_\varepsilon}\|u\|^2_\lambda$, then it is a (PS)-sequence of the free functional $J_\varepsilon$ on $H^s_0(\Omega)$.
\end{lemma}

The above lemma means that the Nehari manifold is a natural constraint for $J_\varepsilon$, whose proof is standard and is omitted (see \cite{Ambrosetti} for details). For simplicity, we write $J_\varepsilon$ instead of $J_\varepsilon|_{\mathcal N_\varepsilon}$.

\section{Two Limit problems}

In this section, we track $\lim\limits_{\varepsilon\to0}\vartheta_\varepsilon$, for which we consider two limit problems.

The first limit problem is
\begin{equation}\label{limit1}
    \begin{cases}
    (-\Delta)^s u=\bigg(\displaystyle\int_{\mathbb R^N}\frac{u^{2^*_{\sigma,s}}(y)}{|x-y|^\sigma}\text{d}y\bigg)u^{2^*_{\sigma,s}-1} \quad &\text{in }\mathbb R^N\\
    u>0 \quad &\text{in }\mathbb R^N
    \end{cases}.
\end{equation}
For \eqref{limit1}, we define the corresponding functional $J_*: D^{s,2}(\mathbb R^N)\to\mathbb R$ by
\begin{equation*}
    J_*(u):=\frac{1}{2}\int_{\mathbb R^{2N}}\frac{|u(x)-u(y)|^2}{|x-y|^{N+2s}}\text{d}x\text{d}y-\frac{1}{2\cdot 2^*_{\sigma,s}}\int_{\mathbb R^{2N}}\frac{(u^+(x))^{2^*_{\sigma,s}}(u^+(y))^{2^*_{\sigma,s}}}{|x-y|^\sigma}\text{d}x\text{d}y,
\end{equation*}
the Nehari manifold associated to $J_*$ by
\begin{equation*}
    \mathcal N_*:=\left\{u\in D^{s,2}(\mathbb R^N)\setminus\{0\}:\left<J'_*(u),u\right>=0\right\},
\end{equation*}
and the least energy by
\begin{equation*}
   \vartheta_*:=\inf_{u\in\mathcal N_*}J_*(u).
\end{equation*}

In the following result, we compute the relation between $\vartheta_*$ and $S^s_{H,L}$ given in \eqref{S1}.

\begin{lemma}\label{limit-energy}
	There holds
	\begin{equation*}
	\vartheta_*=\left(\frac{2^*_{\sigma,s}-1}{2\cdot 2^*_{\sigma,s}}\right)(S^s_{H,L})^{\frac{2^*_{\sigma,s}}{2^*_{\sigma,s}-1}},
	\end{equation*}
	and $\vartheta_*$ is achieved only by functions $(S^s_{H,L})^{\frac{1}{2\cdot2^*_{\sigma,s}-2}}U_{R,a}$ with $U_{R,a}$ defined in \eqref{U}.
\end{lemma}

\begin{proof}[\bf Proof]
	For $u\in D^{s,2}(\mathbb R^N)\setminus\{0\}$, we have
	\begin{equation*}
	\max_{t>0}J_*(tu)=\left(\frac{2^*_{\sigma,s}-1}{2\cdot 2^*_{\sigma,s}}\right)\left(\frac{\displaystyle\int_{\mathbb R^{2N}}\frac{|u(x)-u(y)|^2}{|x-y|^{N+2s}}\text{d}x\text{d}y}{\bigg(\displaystyle\int_{\mathbb R^{2N}}\frac{(u^+(x))^{2^*_{\sigma,s}}(u^+(y))^{2^*_{\sigma,s}}}{|x-y|^\sigma}\text{d}x\text{d}y\bigg)^{\frac{1}{2^*_{\sigma,s}}}}\right)^{\frac{2^*_{\sigma,s}}{2^*_{\sigma,s}-1}}.
	\end{equation*}
	Therefore, $\vartheta_*=\inf_{u\in D^{s,2}(\mathbb R^N)\setminus\{0\}}\max_{t\ge0}J_*(tu)\ge \left(\frac{2^*_{\sigma,s}-1}{2\cdot 2^*_{\sigma,s}}\right)(S^s_{H,L})^{\frac{2^*_{\sigma,s}}{2^*_{\sigma,s}-1}}$.
	
	Additionally, notice that $\tilde u:=(S^s_{H,L})^{\frac{1}{2\cdot2^*_{\sigma,s}-2}}U_{R,a}$ satisfies \eqref{limit1}, we get that $\tilde u\in\mathcal N_*$ and
	\begin{equation*}
	\vartheta_*\le J_*(\tilde u)=\left(\frac{1}{2}-\frac{1}{2\cdot 2^*_{\sigma,s}}\right) \int_{\mathbb R^{2N}}\frac{|\tilde u(x)-\tilde u(y)|^2}{|x-y|^{N+2s}}\text{d}x\text{d}y=\left(\frac{2^*_{\sigma,s}-1}{2\cdot 2^*_{\sigma,s}}\right)(S^s_{H,L})^{\frac{2^*_{\sigma,s}}{2^*_{\sigma,s}-1}}.
	\end{equation*}
\end{proof}

We want to use the minimizers of $\vartheta_*$ to construct the approximating sequences for $\vartheta_\varepsilon$. Since $\Omega$ is a bounded domain in $\mathbb R^N$ with Lipschitz boundary, we choose $r>0$ small enough such that 
\begin{equation*}
    \Omega^+_r:=\{x\in\mathbb R^N:\text{dist}(x,\Omega)\le r\}
\end{equation*}
and  
\begin{equation*}
    \Omega^-_r:=\{x\in\Omega:\text{dist}(x,\partial\Omega)\ge r\}
\end{equation*}
are homotopically equivalent to $\Omega$. For $R>1$ and $x_0\in\Omega^-_r$, define
\begin{equation}\label{u}
    u_{R,x_0}(\cdot):=R^{\frac{N-2s}{2}}U_{1,0}(R(\cdot-x_0))\phi_{r,x_0}(\cdot)=U_{R,x_0}(\cdot)\phi_{r,x_0}(\cdot),
\end{equation}
where $U_{1,0}$ is the standard bubble function defined in \eqref{U} with $R=1$ and $a=0$, and $\phi_{r,x_0}$ is a cut-off function defined by
\begin{equation}\label{cutoff}
   \phi_{r,x_0}(x):=\begin{cases}
   1\quad&\text{if }|x-x_0|<\frac{r}{2}\\
   (0,1)\quad&\text{if }\frac{r}{2}\le|x-x_0|\le r\\
   0\quad&\text{if }|x-x_0|>r
   \end{cases}.
\end{equation}
It follows from $x_0\in\Omega^-_r$ that $u_{R,x_0}\in H^s_0(\Omega)$.

Noticing $N>4s$ and referring to Proposition 21 in Servadei et al \cite{Servadei-Valdinoci}, we can verify that
\begin{align}
    \label{1}
    &\int_\Upsilon \frac{|u_{R,x_0}(x)-u_{R,x_0}(y)|^2}{|x-y|^{N+2s}}\text{d}x\text{d}y\le\int_{\mathbb R^{2N}}\frac{|U_{1,0}(x)-U_{1,0}(y)|^2}{|x-y|^{N+2s}}\text{d}x\text{d}y+o_R(1),\\
    \label{2}
    &\int_{\Omega\times\Omega}\frac{|u_{R,x_0}(x)|^{2^*_{\sigma,s}}|u_{R,x_0}(y)|^{2^*_{\sigma,s}}}{|x-y|^\sigma}\text{d}x\text{d}y=\int_{\mathbb R^{2N}}\frac{|U_{1,0}(x)|^{2^*_{\sigma,s}}|U_{1,0}(y)|^{2^*_{\sigma,s}}}{|x-y|^\sigma}\text{d}x\text{d}y+o_R(1),\\
    \label{3}
    &\int_\Omega|u_{R,x_0}(x)|^2\text{d}x=\frac{1}{R^{2s}}\int_{\mathbb R^N}|U_{1,0}(x)|^2\text{d}x+o_R(1)=o_R(1),
\end{align}
where $o_R(1)$ denotes the quantities that tend to 0 as $R\to+\infty$.

\begin{lemma}\label{<}
	 $\limsup\limits_{\varepsilon\to0}\vartheta_\varepsilon\le\vartheta_*$.
\end{lemma}

\begin{proof}[\bf Proof]
	For any $\varepsilon>0$, there is a unique $t_\varepsilon(u_{R,x_0})>0$ satisfying $t_\varepsilon(u_{R,x_0})u_{R,x_0}\in\mathcal N_\varepsilon$. Thus
	\begin{equation*}
	\|t_\varepsilon(u_{R,x_0})u_{R,x_0}\|^2_\lambda=\int_{\Omega\times\Omega}\frac{|t_\varepsilon(u_{R,x_0})u_{R,x_0}(x)|^{p_\varepsilon}|t_\varepsilon(u_{R,x_0})u_{R,x_0}(y)|^{p_\varepsilon}}{|x-y|^\sigma}\text{d}x\text{d}y,
	\end{equation*}
	which implies
	\begin{equation*}
	t^{2p_\varepsilon-2}_\varepsilon(u_{R,x_0})=\frac{\|u_{R,x_0}\|^2_\lambda}{\displaystyle\int_{\Omega\times\Omega}\frac{|u_{R,x_0}(x)|^{p_\varepsilon}|u_{R,x_0}(y)|^{p_\varepsilon}}{|x-y|^\sigma}\text{d}x\text{d}y}.
	\end{equation*}
	In virtue of \eqref{1}, \eqref{2} and \eqref{3}, we get 
	\begin{equation}\label{tofu}
	\begin{aligned}
	\lim_{\varepsilon\to0}t_\varepsilon(u_{R,x_0})=&\left(\frac{\|u_{R,x_0}\|^2_\lambda}{\displaystyle\int_{\Omega\times\Omega}\frac{|u_{R,x_0}(x)|^{2^*_{\sigma,s}}|u_{R,x_0}(y)|^{2^*_{\sigma,s}}}{|x-y|^\sigma}\text{d}x\text{d}y}\right)^{\frac{1}{2\cdot2^*_{\sigma,s}-2}}\\
	\le&\left(\frac{\|U_{1,0}\|^2_{s,2}+o_R(1)}{\displaystyle\int_{\mathbb R^{2N}}\frac{|U_{1,0}(x)|^{2^*_{\sigma,s}}|U_{1,0}(y)|^{2^*_{\sigma,s}}}{|x-y|^\sigma}\text{d}x\text{d}y+o_R(1)}\right)^{\frac{1}{2\cdot2^*_{\sigma,s}-2}}\\
	=&(S^s_{H,L})^{\frac{1}{2\cdot2^*_{\sigma,s}-2}}+o_R(1),
	\end{aligned}
    \end{equation}
	and then
	\begin{equation*}
	J_\varepsilon(t_\varepsilon(u_{R,x_0})u_{R,x_0})=\frac{p_\varepsilon-1}{2p_\varepsilon}\|t_\varepsilon(u_{R,x_0})u_{R,x_0}\|^2_\lambda=\frac{p_\varepsilon-1}{2p_\varepsilon}t^2_\varepsilon(u_{R,x_0})\|U_{1,0}\|^2_{s,2}+o_R(1).
	\end{equation*}
	Therefore,
	\begin{equation*}
	\lim_{\varepsilon\to0}J_\varepsilon(t_\varepsilon(u_{R,x_0})u_{R,x_0})\le\frac{2^*_{\sigma,s}-1}{2\cdot2^*_{\sigma,s}}(S^s_{H,L})^{\frac{1}{2^*_{\sigma,s}-1}}\|U_{1,0}\|^2_{s,2}+o_R(1).
	\end{equation*}
	For any $\delta>0$, we can choose $R$ large enough such that $o_R(1)<\delta$. By Lemma \ref{limit-energy}, we obtain that
	\begin{equation*}
	\limsup_{\varepsilon\to0}\vartheta_\varepsilon\le\lim_{\varepsilon\to0}J_\varepsilon(t_\varepsilon(u_{R,x_0})u_{R,x_0})<\frac{2^*_{\sigma,s}-1}{2\cdot2^*_{\sigma,s}}(S^s_{H,L})^{\frac{1}{2^*_{\sigma,s}-1}}\|U_{1,0}\|^2_{s,2}+\delta=\vartheta_*+\delta.
	\end{equation*}
	This completes the proof by letting $\delta\to0$.	
\end{proof}

\begin{remark}\label{bound}
	The groundstates $u_\varepsilon$ are also bounded  uniformly in $\varepsilon$. Indeed, by Lemma \ref{groundstate},
	\begin{equation*}
	\|u_\varepsilon\|^2_\lambda=\frac{2p_\varepsilon}{p_\varepsilon-1}J_\varepsilon(u_\varepsilon)=\frac{2p_\varepsilon}{p_\varepsilon-1}\vartheta_\varepsilon.
	\end{equation*}
\end{remark}

Now we introduce the second limit problem which acts as the mediator between the problem \eqref{target} and the first limit problem \eqref{limit1}. In particular, it will play an important role in computing $\lim\limits_{\varepsilon\to0}\vartheta_\varepsilon$. Consider
\begin{equation}\label{limit2}
    \begin{cases}
    (-\Delta)^s u+\lambda u=\bigg(\displaystyle\int_\Omega\frac{u^{2^*_{\sigma,s}}(y)}{|x-y|^\sigma}\text{d}y\bigg)u^{2^*_{\sigma,s}-1} \quad &\text{in }\Omega\\
    u=0 \quad &\text{in }\mathbb R^N\setminus\Omega\\
    u>0 \quad &\text{in }\Omega
    \end{cases}.
\end{equation}
The existence, nonexistence and regularity results of weak solutions to \eqref{limit2} have been studied in Mukherjee et al \cite{Mukherjee-Sreenadh}. As usual, we define the energy functional $J^\Omega_*:H^s_0(\Omega)\to\mathbb R$ for \eqref{limit2} by
\begin{equation*}
    J^\Omega_*(u):=\frac{1}{2}\|u\|^2_\lambda-\frac{1}{2\cdot2^*_{\sigma,s}}\int_{\Omega\times\Omega}\frac{(u^+(x))^{2^*_{\sigma,s}}(u^+(y))^{2^*_{\sigma,s}}}{|x-y|^\sigma}\text{d}x\text{d}y,
\end{equation*}
the associated Nehari manifold by
\begin{equation*}
\mathcal N^\Omega_*:=\left\{u\in H^s_0(\Omega)\setminus\{0\}:\left<(J^\Omega_*)'(u),u\right>=0\right\},
\end{equation*}
and the least energy by
\begin{equation*}
\vartheta^\Omega_*:=\inf_{u\in\mathcal N^\Omega_*}J^\Omega_*(u).
\end{equation*}

\begin{lemma}\label{not}
	$\vartheta^\Omega_*=\vartheta_*$ and $\vartheta^\Omega_*$ is not achieved.
\end{lemma}

\begin{proof}[\bf Proof]
	We first show $\vartheta^\Omega_*=\vartheta_*$. For one thing, for any $u\in \mathcal N^\Omega_*$, we extend $u$ to zero outside $\Omega$, and then there is a unique $t_*(u)\in(0,1)$ such that $t_*(u)u\in\mathcal N_*$. Hence
	\begin{equation*}
	\vartheta_*\le J_*(t_*(u)u)=\left(\frac{1}{2}-\frac{1}{2\cdot2^*_{\sigma,s}}\right)\|t_*(u)u\|^2_{s,2}<\left(\frac{1}{2}-\frac{1}{2\cdot2^*_{\sigma,s}}\right)\|u\|^2_\lambda,
	\end{equation*}
	which implies $\vartheta_*\le \vartheta^\Omega_*$. For another, for all $x_0\in\Omega^-_r$ and $R>1$, we take $u_{R,x_0}$ defined in \eqref{u} and a unique $t^\Omega_*(u_{R,x_0})>0$ satisfying $t^\Omega_*(u_{R,x_0})u_{R,x_0}\in\mathcal N^\Omega_*$. Proceeding as the proof of Lemma \ref{<}, we obtain that for any $\delta>0$, there exists $R>1$ such that
	\begin{equation*}
	\vartheta^\Omega_*\le J^\Omega_*(t^\Omega_*(u_{R,x_0})u_{R,x_0})<\vartheta_*+\delta.
	\end{equation*}
	We get $\vartheta^\Omega_*\le\vartheta_*$ by means of the arbitrariness of $\delta$. Hence $\vartheta^\Omega_*=\vartheta_*$.
	
	Next we show that $\vartheta^\Omega_*$ can not be achieved. Indeed, suppose by contradiction that $v\in \mathcal N^\Omega_*$ satisfies $J^\Omega_*(v)=\vartheta^\Omega_*$. 
	We extend $v$ to zero outside $\Omega$. There is a unique $t_*(v)>0$ satisfying $t_*(v)v\in\mathcal N_*$. Thus
	\begin{align*}
	t^2_*(v)\|v\|^2_{s,2}=&t^{2\cdot 2^*{\sigma,s}}_*(v)\int_{\mathbb R^{2N}}\frac{(v^+(x))^{2^*_{\sigma,s}}(v^+(y))^{2^*_{\sigma,s}}}{|x-y|^\sigma}\text{d}x\text{d}y\\
	=&t^{2\cdot 2^*{\sigma,s}}_*(v)\int_{\Omega\times\Omega}\frac{(v^+(x))^{2^*_{\sigma,s}}(v^+(y))^{2^*_{\sigma,s}}}{|x-y|^\sigma}\text{d}x\text{d}y\\
	=&t^{2\cdot 2^*{\sigma,s}}_*(v)\|v\|^2_\lambda,
	\end{align*}
	which together with $\lambda>0$ implies that 
	$$t_*(v)=\left(\frac{\|v\|^2_{s,2}}{\|v\|^2_\lambda}\right)^{\frac{1}{2\cdot 2^*{\sigma,s}-2}}<1.$$
	Hence
	\begin{align*}
	\vartheta_*\le J_*(t_*(v)v)=&\left(\frac{1}{2}-\frac{1}{2\cdot2^*_{\sigma,s}}\right)t^{2\cdot 2^*{\sigma,s}}_*(v)\int_{\mathbb R^{2N}}\frac{(v^+(x))^{2^*_{\sigma,s}}(v^+(y))^{2^*_{\sigma,s}}}{|x-y|^\sigma}\text{d}x\text{d}y\\
	=&\left(\frac{1}{2}-\frac{1}{2\cdot2^*_{\sigma,s}}\right)t^{2\cdot 2^*{\sigma,s}}_*(v)\int_{\Omega\times\Omega}\frac{(v^+(x))^{2^*_{\sigma,s}}(v^+(y))^{2^*_{\sigma,s}}}{|x-y|^\sigma}\text{d}x\text{d}y\\
	<&\left(\frac{1}{2}-\frac{1}{2\cdot2^*_{\sigma,s}}\right)\int_{\Omega\times\Omega}\frac{(v^+(x))^{2^*_{\sigma,s}}(v^+(y))^{2^*_{\sigma,s}}}{|x-y|^\sigma}\text{d}x\text{d}y=\vartheta^\Omega_*=\vartheta_*,
	\end{align*}
	which is a contradiction.
	
	If $\lambda=0$, then it follows from $v\in \mathcal N^\Omega_*$ that $t_*(v)=1$ and then $v\in\mathcal N_*$. Moreover, $$J_*(v)=J^\Omega_*(v)=\vartheta^\Omega_*=\vartheta_*=\inf_{u\in\mathcal N_*}J_*(u),$$
	which implies that $v=(S^s_{H,L})^{\frac{1}{2\cdot2^*_{\sigma,s}-2}}U_{R,a}>0$. But it is impossible by the construction of $v$. 
\end{proof}

Now we prove the main result in this section.

\begin{proposition}\label{=}
	There holds $$\lim\limits_{\varepsilon\to0}\vartheta_\varepsilon=\vartheta_*.$$
\end{proposition}

\begin{proof}[\bf Proof]
	In virtue of Lemma \ref{<}, it is sufficient to show that
	\begin{equation}\label{>}
	\liminf_{\varepsilon\to0}\vartheta_\varepsilon\ge\vartheta_*.
	\end{equation}
	
	By Lemma \ref{groundstate}, we choose $u_\varepsilon\in\mathcal N_\varepsilon$ satisfying $J_\varepsilon(u_\varepsilon)=\vartheta_\varepsilon$. There is a unique $t^\Omega_*(u_\varepsilon)>0$ such that $t^\Omega_*(u_\varepsilon)u_\varepsilon\in\mathcal N^\Omega_*$, namely,
	\begin{equation*}
	\left(t^\Omega_*(u_\varepsilon)\right)^2\|u_\varepsilon\|^2_\lambda=\left(t^\Omega_*(u_\varepsilon)\right)^{2\cdot2^*_{\sigma,s}}\int_{\Omega\times\Omega}\frac{(u_\varepsilon^+(x))^{2^*_{\sigma,s}}(u_\varepsilon^+(y))^{2^*_{\sigma,s}}}{|x-y|^\sigma}\text{d}x\text{d}y.
	\end{equation*}
	Noticing $u_\varepsilon\in\mathcal N_\varepsilon$, we get that
	\begin{equation}\label{t*1}
	\begin{aligned}
	\left(t^\Omega_*(u_\varepsilon)\right)^{2\cdot2^*_{\sigma,s}-2}=&\frac{\|u_\varepsilon\|^2_\lambda}{\displaystyle\int_{\Omega\times\Omega}\frac{(u_\varepsilon^+(x))^{2^*_{\sigma,s}}(u_\varepsilon^+(y))^{2^*_{\sigma,s}}}{|x-y|^\sigma}\text{d}x\text{d}y}\\
	=&\frac{\displaystyle\int_{\Omega\times\Omega}\frac{(u_\varepsilon^+(x))^{p_\varepsilon}(u_\varepsilon^+(y))^{p_\varepsilon}}{|x-y|^\sigma}\text{d}x\text{d}y}{\displaystyle\int_{\Omega\times\Omega}\frac{(u_\varepsilon^+(x))^{2^*_{\sigma,s}}(u_\varepsilon^+(y))^{2^*_{\sigma,s}}}{|x-y|^\sigma}\text{d}x\text{d}y}.
	\end{aligned}
	\end{equation} 
	We claim that 
	\begin{equation}\label{<1}
	\limsup_{\varepsilon\to0}t^\Omega_*(u_\varepsilon)\le 1.
	\end{equation}
	
	It follows from
	\begin{equation*}
	\int_{\Omega\times\Omega}\frac{(u_\varepsilon^+(x))^{p_\varepsilon}(u_\varepsilon^+(y))^{p_\varepsilon}}{|x-y|^\sigma}\text{d}x\text{d}y=\int_{\Omega\times\Omega}\frac{(u_\varepsilon^+(x))^{p_\varepsilon}(u_\varepsilon^+(y))^{p_\varepsilon}}{|x-y|^{\sigma\cdot\frac{p_\varepsilon}{2^*_{\sigma,s}}}}\cdot\frac{1}{|x-y|^{\sigma\cdot\frac{\varepsilon}{2^*_{\sigma,s}}}}\text{d}x\text{d}y
	\end{equation*}
	and the H\"older inequality that 
	\begin{align*}
	&\int_{\Omega\times\Omega}\frac{(u_\varepsilon^+(x))^{p_\varepsilon}(u_\varepsilon^+(y))^{p_\varepsilon}}{|x-y|^\sigma}\text{d}x\text{d}y\\
	\le&\bigg(\int_{\Omega\times\Omega}\frac{(u_\varepsilon^+(x))^{2^*_{\sigma,s}}(u_\varepsilon^+(y))^{2^*_{\sigma,s}}}{|x-y|^\sigma}\text{d}x\text{d}y\bigg)^{\frac{p_\varepsilon}{2^*_{\sigma,s}}}\cdot\left(\int_{\Omega\times\Omega}\frac{1}{|x-y|^\sigma}\text{d}x\text{d}y\right)^{\frac{\varepsilon}{2^*_{\sigma,s}}}.
	\end{align*}
	We take a change of variables $\xi=x-y$, $\eta=x+y$. Then for $\rho=\rho(\Omega)>0$ large enough, we have
	\begin{equation*}
	\int_{\Omega\times\Omega}\frac{1}{|x-y|^\sigma}\text{d}x\text{d}y\le\frac{1}{2}\int_{B_\rho(0)\times B_\rho(0)}\frac{1}{|\xi|^\sigma}\text{d}\xi\text{d}\eta\le C_\rho\int_{B_\rho(0)}\frac{1}{|\xi|^\sigma}\text{d}\xi=C_\rho,
	\end{equation*}
	where we have used the assumption $\sigma\in(0,N)$. Since $\rho$ depends only on $\Omega$, we obtain that
	\begin{equation}\label{p}
	\int_{\Omega\times\Omega}\frac{(u_\varepsilon^+(x))^{p_\varepsilon}(u_\varepsilon^+(y))^{p_\varepsilon}}{|x-y|^\sigma}\text{d}x\text{d}y\le\bigg(\int_{\Omega\times\Omega}\frac{(u_\varepsilon^+(x))^{2^*_{\sigma,s}}(u_\varepsilon^+(y))^{2^*_{\sigma,s}}}{|x-y|^\sigma}\text{d}x\text{d}y\bigg)^{\frac{p_\varepsilon}{2^*_{\sigma,s}}}\cdot C_\Omega^{\frac{\varepsilon}{2^*_{\sigma,s}}}.
	\end{equation}
	By inserting \eqref{p} into \eqref{t*1}, we get that
	\begin{equation}\label{t*2}
	\begin{aligned}
	\left(t^\Omega_*(u_\varepsilon)\right)^{2\cdot2^*_{\sigma,s}-2}
	\le&\bigg(\displaystyle\int_{\Omega\times\Omega}\frac{(u_\varepsilon^+(x))^{2^*_{\sigma,s}}(u_\varepsilon^+(y))^{2^*_{\sigma,s}}}{|x-y|^\sigma}\text{d}x\text{d}y\bigg)^{\frac{-\varepsilon}{2^*_{\sigma,s}}}\cdot C_\Omega^{\frac{\varepsilon}{2^*_{\sigma,s}}}\\
	=&\left(\frac{C_\Omega}{\int_{\Omega\times\Omega}\frac{(u_\varepsilon^+(x))^{2^*_{\sigma,s}}(u_\varepsilon^+(y))^{2^*_{\sigma,s}}}{|x-y|^\sigma}\text{d}x\text{d}y}\right)^{\frac{\varepsilon}{2^*_{\sigma,s}}}.
	\end{aligned}
	\end{equation}
	Due to Remark \ref{bound} and \eqref{SHL}, we see that $\int_{\Omega\times\Omega}\frac{(u_\varepsilon^+(x))^{2^*_{\sigma,s}}(u_\varepsilon^+(y))^{2^*_{\sigma,s}}}{|x-y|^\sigma}\text{d}x\text{d}y$ is bounded uniformly in $\varepsilon$. Hence we deduce \eqref{<1} from \eqref{t*2}. 
	
	Concequently, 
	\begin{align*}
	\vartheta_*&=\vartheta^\Omega_*\le J^\Omega_*(t^\Omega_*(u_\varepsilon)u_\varepsilon)=\left(\frac{1}{2}-\frac{1}{2\cdot2^*_{\sigma,s}}\right)\left(t^\Omega_*(u_\varepsilon)\right)^2\|u_\varepsilon\|^2_\lambda\\
	&=\left(t^\Omega_*(u_\varepsilon)\right)^2\frac{\frac{1}{2}-\frac{1}{2\cdot2^*_{\sigma,s}}}{\frac{1}{2}-\frac{1}{2p_\varepsilon}}\left(\frac{1}{2}-\frac{1}{2p_\varepsilon}\right)\|u_\varepsilon\|^2_\lambda\\
	&\le(1+o_\varepsilon(1))\vartheta_\varepsilon,
	\end{align*}
	where $o_\varepsilon(1)\to 0$ as $\varepsilon\to0$ and $\vartheta_\varepsilon$ is bounded by Lemma \ref{<}. Thus \eqref{>} is showed.
\end{proof}

In the end of this section, we state a technical lemma from He and R$\breve a$dulescu \cite{He}, and prove a nonlocal splitting lemma which gives a complete description for the functional $J^\Omega_*$. This nonlocal splitting lemma is a variant of the classical one contained in Struwe \cite{Struwe}.

\begin{lemma}\label{He}
	(Lemma 3.1 in \cite{He})
	Let $\{v_n\}_{n\in\mathbb N}$ be a (PS)$_c$-sequence for the functional $J_*$ with $v_n\rightharpoonup0$ and $v_n\nrightarrow0$ in $D^{s,2}(\mathbb R^N)$ as $n\to\infty$. Then there exists a sequence $\{R_n\}_{n\in\mathbb N}\subset\mathbb R^+$, a point sequence $\{x_n\}_{n\in\mathbb N}\subset\mathbb R^N$ and a nontrival solution $v_0\in D^{s,2}(\mathbb R^N)$ of \eqref{limit1} such that, up to a subsequence of $\{v_n\}_{n\in\mathbb N}$, we have that
	$$\tilde v_n(x)=v_n(x)-R_n^{\frac{N-2s}{2}}v_0(R_n(x-x_n))+o_n(1)$$
	is a (PS)$_{c-J_*(v_0)}$-sequence for $J_*$, where $o_n(1)\to0$ as $n\to\infty$.
\end{lemma}

The above property on (PS)-sequence to the first limit problem \eqref{limit1} is very important in proving the following nonlocal splitting lemma for the fractional critical Choquard problem \eqref{limit2}.

\begin{lemma}\label{split}
	Let $\{v_n\}_{n\in\mathbb N}$ be a (PS)-sequence of $J^\Omega_*$ in $H^s_0(\Omega)$. Then there exist $k\in\mathbb N$, a point sequence $\{x^j_n\}_{n\in\mathbb N}\subset\Omega$, a radius sequence $\{R^j_n\}_{n\in\mathbb N}\subset\mathbb R^+$, a solution $v\in H^s_0(\Omega)$ of \eqref{limit2}, and nontrivial solutions $v^j\in D^{s,2}(\mathbb R^N)$ to \eqref{limit1}, where $j=1,2,\cdots,k$, such that a subsequence of $\{v_n\}_{n\in\mathbb N}$, denoted also by $\{v_n\}_{n\in\mathbb N}$, satisfies 
	\begin{equation*}
	\Big\|v_n-v-\sum_{j=1}^k v^j_{R^j_n,x^j_n}\Big\|_{s,2}\to0\quad\text{as }n\to\infty,
	\end{equation*}	 
	where 
    \begin{equation*}
    v^j_{R^j_n,x^j_n}(x):=(R^j_n)^{\frac{N-2s}{2}}v^j(R^j_n(x-x^j_n)), \quad j=1,2,\cdots,k.
    \end{equation*}
    Moreover,
    \begin{equation}
    J^\Omega_*(v_n)\to J^\Omega_*(v)+\sum_{j=1}^k
    J_*(v^j)\quad\text{as }n\to\infty.    
    \end{equation}
\end{lemma}

\begin{proof}[\bf Proof]
	{\bf Step 1.}	
	The (PS)-sequence $\{v_n\}_{n\in\mathbb N}$ of $J^\Omega_*$ is bounded in $H^s_0(\Omega).$ Thus up to a subsequence, we may assume that $v_n\rightharpoonup v$ in $H^s_0(\Omega)$ as $n\to\infty$. Moreover, $v$ solves \eqref{limit2}. Set $v^1_n:=v_n-v$. Then by Br\'ezis-Lieb Lemma (see Lemma 1.32 and Remark 1.33 in Willem \cite{Willem}), we get
	\begin{align*}
	\|v_n\|^2_{s,2}-\|v^1_n\|^2_{s,2}&\to\|v\|^2_{s,2}\quad\text{as }n\to\infty,\\
	|v_n|^2_2-|v^1_n|^2_2&\to|v|^2_2\quad\text{as }n\to\infty,\\
	(v^+_n)^{2^*{\sigma,s}}-((v^1_n)^+)^{2^*{\sigma,s}}&\to(v^+)^{2^*{\sigma,s}}\quad\text{in } L^{\frac{2N}{2N-\sigma}}(\Omega) \ \text{ as }n\to\infty.
	\end{align*}
	By Lemma 2.1 in Liu et al \cite{Liu-DCDS}, we have
	\begin{align*}
	|x|^{-\sigma}*(v^+_n)^{2^*{\sigma,s}}-|x|^{-\sigma}*((v^1_n)^+)^{2^*{\sigma,s}}\to|x|^{-\sigma}*(v^+)^{2^*{\sigma,s}}\quad\text{in } L^{\frac{2N}{\sigma}}(\Omega) \ \text{ as }n\to\infty.
	\end{align*}
	Proceeding as the arguments of Lemma 2.5 in Liu et al \cite{Liu-DCDS} with slight amendment, we obtain
	\begin{equation}\label{AppendixA3}
	\begin{aligned}
	&\int_{\Omega\times\Omega}\frac{(v^+_n(x))^{2^*{\sigma,s}}(v^+_n(y))^{2^*{\sigma,s}}}{|x-y|^\sigma}\text{d}x\text{d}y-\int_{\Omega\times\Omega}\frac{((v^1_n)^+(x))^{2^*{\sigma,s}}((v^1_n)^+(y))^{2^*{\sigma,s}}}{|x-y|^\sigma}\text{d}x\text{d}y\\
	\to&\int_{\Omega\times\Omega}\frac{(v^+(x))^{2^*{\sigma,s}}(v^+(y))^{2^*{\sigma,s}}}{|x-y|^\sigma}\text{d}x\text{d}y\quad \text{as }n\to\infty,
	\end{aligned}
	\end{equation}
	and 
	\begin{equation}\label{AppendixA9}
	\begin{aligned}
	&\int_{\Omega\times\Omega}\frac{(v^+_n(y))^{2^*{\sigma,s}}(v^+_n(x))^{2^*{\sigma,s}-1}\psi(x)}{|x-y|^\sigma}\text{d}x\text{d}y\\
	&-\int_{\Omega\times\Omega}\frac{((v^1_n)^+(y))^{2^*{\sigma,s}}((v^1_n)^+(x))^{2^*{\sigma,s}-1}\psi(x)}{|x-y|^\sigma}\text{d}x\text{d}y\\
	\to&\int_{\Omega\times\Omega}\frac{(v^+(y))^{2^*{\sigma,s}}(v^+(x))^{2^*{\sigma,s}-1}\psi(x)}{|x-y|^\sigma}\text{d}x\text{d}y \  \text{ as }n\to\infty\text{ uniformly in } \psi\in H^s_0(\Omega).
	\end{aligned}
	\end{equation}
	Indeed, the main differences with the proof of Lemma 2.5 in Liu et al \cite{Liu-DCDS} are
	\begin{align*}
	&\int_{\Omega}\left(|x|^{-\sigma}*(v^+)^{2^*{\sigma,s}}\right)\left((v^1_n)^+(x)\right)^{2^*{\sigma,s}-1}|\psi(x)|\text{d}x\\
	\le&\left(\int_{\Omega}\left(|x|^{-\sigma}*|v|^{2^*{\sigma,s}}\right)^{\frac{2N}{N+2s}}|v^1_n(x)|^{(2^*{\sigma,s}-1)\frac{2N}{N+2s}}\text{d}x\right)^{\frac{N+2s}{2N}}|\psi|_{\frac{2N}{N-2s}}\\
	\le&o_n(1)\|\psi\|_\lambda,
	\end{align*}
	and 
	\begin{align*}
	&\int_{\Omega}\left(|x|^{-\sigma}*\left((v^1_n)^+\right)^{2^*{\sigma,s}}\right)\left(v^+(x)\right)^{2^*{\sigma,s}-1}|\psi(x)|\text{d}x\\
	\le&\left(\int_{\Omega}\left(|x|^{-\sigma}*|v^1_n|^{2^*{\sigma,s}}\right)^{\frac{2N}{N+2s}}|v(x)|^{(2^*{\sigma,s}-1)\frac{2N}{N+2s}}\text{d}x\right)^{\frac{N+2s}{2N}}|\psi|_{\frac{2N}{N-2s}}\\
	\le&o_n(1)\|\psi\|_\lambda, 
	\end{align*}
	where the former is ensured by $(|x|^{-\sigma}*|v|^{2^*{\sigma,s}})^{\frac{2N}{N+2s}}\in L^{\frac{N+2s}{\sigma}}(\Omega)$ and $|v^1_n|^{(2^*{\sigma,s}-1)\frac{2N}{N+2s}}\rightharpoonup 0$ in $L^{\frac{N+2s}{N+2s-\sigma}}(\Omega)$ as $n\to\infty$, and the latter is due to $|v|^{(2^*{\sigma,s}-1)\frac{2N}{N+2s}}\in L^{\frac{N+2s}{N+2s-\sigma}}(\Omega)$ and $|x|^{-\sigma}*|v^1_n|^{2^*{\sigma,s}})^{\frac{2N}{N+2s}}\rightharpoonup 0$ in $L^{\frac{N+2s}{\sigma}}(\Omega)$ as $n\to\infty$.
	Therefore, as $n\to\infty$,
	\begin{align*}
	J^\Omega_*(v_n)-J^\Omega_*(v^1_n)&\to J^\Omega_*(v)\\
	(J^\Omega_*)'(v_n)-(J^\Omega_*)'(v^1_n)&\to (J^\Omega_*)'(v) \quad \text{in } H^{-s}_0(\Omega).
	\end{align*}
	Since $(J^\Omega_*)'(v_n)\to0$ in $H^{-s}_0(\Omega)$ as $n\to\infty$ and $(J^\Omega_*)'(v)=0$, we get
	\begin{align*}
	(J^\Omega_*)'(v^1_n)\to0 \quad \text{in } H^{-s}_0(\Omega) \text{ as } n\to\infty.  
	\end{align*}
	Moreover, it follows from $v^1_n\to0$ in $L^2(\Omega)$ as $n\to\infty$ that
	\begin{equation}\label{AppendixA5}
	\begin{aligned}
	J_*(v^1_n)&=J^\Omega_*(v^1_n)+o_n(1)=J^\Omega_*(v_n)-J^\Omega_*(v)+o_n(1)\\
	J'_*(v^1_n)&=(J^\Omega_*)'(v^1_n)+o_n(1)=o_n(1).
	\end{aligned}
	\end{equation}
	Thus $\{v_n^1\}_{n\in\mathbb N}$ is a (PS)-sequence of $J_*$.
	
	{\bf Step 2.} If $v_n^1\to0$ in $D^{s,2}(\mathbb R^N)$ as $n\to\infty$, then the proof is completed with $k=0$. If $v_n^1\nrightarrow0$ in $D^{s,2}(\mathbb R^N)$ as $n\to\infty$, then it follows from Lemma \ref{He} that there exist $\{R_n^1\}_{n\in\mathbb N}\subset\mathbb R^+$, $\{x_n^1\}_{n\in\mathbb N}\subset\mathbb R^N$ and a nontrival solution $v^1\in D^{s,2}(\mathbb R^N)$ of \eqref{limit1} such that 
	$$v_n^2(x):=v_n^1(x)-(R_n^1)^{\frac{N-2s}{2}}v^1(R_n^1(x-x_n^1))+o_n(1)$$
	is a (PS)-sequence for $J_*$. 	

	Define 
	\begin{align*}
	&\tilde v^1_n(x):=(R^1_n)^{\frac{2s-N}{2}}v^1_n\Big(\frac{x}{R^1_n}+x^1_n\Big),\\
	&\tilde v^2_n(x):=(R^1_n)^{\frac{2s-N}{2}}v^2_n\Big(\frac{x}{R^1_n}+x^1_n\Big).
	\end{align*}
	Obviously, $\tilde v^2_n(x)=\tilde v^1_n(x)-v^1(x)+o_n(1)$. Meanwhile, $\tilde v_n^2\rightharpoonup0$ in $D^{s,2}(\mathbb R^N)$ due to the proof of Lemma 3.1 in \cite{He}. Similarly as before,
	\begin{align*}
	\|v^2_n\|^2_{s,2}=\|\tilde v^2_n\|^2_{s,2}=\|\tilde v_n^1\|^2_{s,2}-\|v^1\|^2_{s,2}+o_n(1)=\|v_n\|^2_{s,2}-\|v\|^2_{s,2}-\|v^1\|^2_{s,2}+o_n(1).
	\end{align*}
	Since $J'_*(v^1)=0$ and $J'_*(\tilde v^1_n)=o_n(1)$, we have $\|J'_*(v^2_n)\|=\|J'_*(\tilde v^2_n)\|=o_n(1)$. In addition, by \eqref{AppendixA5}, we get
	\begin{align*}
	J_*(v^2_n)=J_*(\tilde v^2_n)=J_*(\tilde v^1_n)-J_*(v^1)+o_n(1)=J^\Omega_*(v_n)-J^\Omega_*(v)-J_*(v^1)+o_n(1).
	\end{align*}
	
	If $\tilde v_n^2\to0$ in $D^{s,2}(\mathbb R^N)$ as $n\to\infty$, then $v_n^2\to0$ in $D^{s,2}(\mathbb R^N)$ as $n\to\infty$, and the proof is completed with $k=1$. If $\tilde v_n^2\nrightarrow0$ in $D^{s,2}(\mathbb R^N)$ as $n\to\infty$, then $v_n^2\nrightarrow0$ in $D^{s,2}(\mathbb R^N)$ as $n\to\infty$, and it follows from Lemma \ref{He} that there exist $\{R_n^2\}_{n\in\mathbb N}\subset\mathbb R^+$, $\{x_n^2\}_{n\in\mathbb N}\subset\mathbb R^N$ and a nontrival solution $v^2\in D^{s,2}(\mathbb R^N)$ of \eqref{limit1} such that 
	$$v_n^3(x):=v_n^2(x)-(R_n^2)^{\frac{N-2s}{2}}v^2(R_n^2(x-x_n^2))+o_n(1)$$
	is a (PS)-sequence for $J_*$.
	
	Iterating the above procedure, we construct sequences $v^j$, $x^j_n$ and $R^j_n$ such that
	\begin{align*}
	&v^k_n:=v_n-v-\sum_{j=1}^{k-1}(R^j_n)^{\frac{N-2s}{2}}v^j\left(R^j_n(\cdot-x^j_n)\right)+o_n(1),\\
	&\|v^k_n\|^2_{s,2}=\|v_n\|^2_{s,2}-\|v\|^2_{s,2}-\sum_{j=1}^{k-1}\|v^j\|^2_{s,2}+o_n(1),\\
	&J_*(v^k_n)=J^\Omega_*(v_n)-J^\Omega_*(v)-\sum_{j=1}^{k-1}J_*(v^j)+o_n(1),\\
	&J'_*(v^k_n)=o_n(1),\quad J'_*(v^j)=0,\quad j=1,2,\cdots, k-1.
	\end{align*}
	
	{\bf Step 3.} By means of \eqref{FSI} and \eqref{critical<}, we obtain that any nontrival critical point $u$ of $J_*$ satisfies
	\begin{align*}
	|u|^2_{2^*_s}\le S_{N,s}\|u\|^2_{s,2}=S_{N,s}\int_{\mathbb R^{2N}}\frac{(u^+(x))^{2^*_{\sigma,s}}(u^+(y))^{2^*_{\sigma,s}}}{|x-y|^\sigma}\text{d}x\text{d}y\le S_{N,s} C_{N,\sigma,s}|u|^{2\cdot 2^*_{\sigma,s}}_{2^*_s},
	\end{align*}
	which implies that
	\begin{align}\label{AppendixA8}
	J_*(u)=\bigg(\frac{1}{2}-\frac{1}{2\cdot 2^*_{\sigma,s}}\bigg)\|u\|^2_{s,2}\ge\frac{2^*_{\sigma,s}-1}{2\cdot 2^*_{\sigma,s}}\frac{1}{S_{N,s}}(S^s_{H,L})^{\frac{1}{2^*_{\sigma,s}-1}}>0.
	\end{align}
	
	Theorefore, the above iteration must terminate at some finite index $k$ by \eqref{AppendixA8}. Hence $\|v^k_n\|_{s,2}=\|\tilde v^k_n\|_{s,2}\to0$ as $n\to\infty$. The proof of Lemma \ref{split} is finished.
\end{proof}

In virtue of the above nonlocal splitting lemma and Lemmas \ref{limit-energy}, \ref{not}, we get the following immediate result.

\begin{remark}\label{k=1}
	If there exists a (PS)-sequence for $J^\Omega_*$ at level $\vartheta^\Omega_*$, then $$v=0,\quad k=1,\quad v^1=(S^s_{H,L})^{\frac{1}{2\cdot2^*_{\sigma,s}-2}}U_{R,a}$$ and $v_n-(S^s_{H,L})^{\frac{1}{2\cdot2^*_{\sigma,s}-2}}U_{R_n,a_n}\to0$ in $D^{s,2}(\mathbb R^n)$ as $n\to\infty$.
\end{remark}

\section{Proof of the main result}

By applying Proposition \ref{BCP} to prove Theorem \ref{t}, we need to construct a map from $\Omega^-_r$ to $\mathcal N_\varepsilon$ and a function from $\mathcal N_\varepsilon$ to $\Omega^+_r$. Denote with the same symbol $u$ its trivial extension out of the support of $u$. We introduce the barycenter of a function $u\in D^{s,2}(\mathbb R^N)$ with compact support as $\beta(u):=(\beta^1(u),\beta^2(u),\cdots,\beta^N(u))$, where
\begin{equation}\label{barycenter}
    \beta^i(u):=\frac{\displaystyle\int_{\mathbb R^N}x^i|u|^{2^*_s}\text{d}x}{\displaystyle\int_{\mathbb R^N}|u|^{2^*_s}\text{d}x},\quad i=1,2,\cdots N.
\end{equation}  
This barycenter map allows us to compare the topology of $\Omega$ with the topology of some suitable sublevels of $J_\varepsilon$. Exectly, we can show the following result  according to Remark \ref{k=1}.

\begin{proposition}\label{beta}
	There exist $\delta_0>0$ and $\varepsilon_0=\varepsilon_0(\delta_0)>0$ such that for any $\delta\in(0,\delta_0]$ and for any $\varepsilon\in(0,\varepsilon_0]$, it holds
	$$u\in\mathcal N_\varepsilon\ \text{ and }\ J_\varepsilon(u)<\vartheta_\varepsilon+\delta\Longrightarrow \beta(u)\in\Omega^+_r.$$
\end{proposition}

\begin{proof}[\bf Proof]
	Suppose on the contrary that there exist sequences $\delta_n\to0$, $\varepsilon_n\to0$ and $u_n\in\mathcal N_{\varepsilon_n}$ such that
	\begin{equation}\label{contradiction}
	J_{\varepsilon_n}(u_n)\le \vartheta_{\varepsilon_n}+\delta_n\quad\text{and}\quad\beta(u_n)\notin\Omega^+_r.
	\end{equation}
	It follows from \eqref{contradiction} and Proposition \ref{=} that
	\begin{equation}\label{nto}
	J_{\varepsilon_n}(u_n)\to\vartheta_*\quad\text{as }n\to\infty,
	\end{equation}
	and $\{u_n\}_{n\in\mathbb N}$ is bounded in $H^s_0(\Omega)$. There is a unique $t^\Omega_{*,n}(u_n)>0$ such that $t^\Omega_{*,n}(u_n)u_n\in\mathcal N^\Omega_*$. Set $p_n:=2^*_{\sigma,s}-\varepsilon_n$, we next evaluate
	\begin{align*}
	&J_{\varepsilon_n}(u_n)-J^\Omega_*\left(t^\Omega_{*,n}(u_n)u_n\right)\\
	=&\left(\frac{1}{2}-\frac{1}{2p_n}\right)\|u_n\|^2_\lambda-\left(\frac{1}{2}-\frac{1}{2\cdot 2^*_{\sigma,s}}\right)\left(t^\Omega_{*,n}(u_n)\right)^2\|u_n\|^2_\lambda\\
	=&\left(\frac{1}{2}-\frac{1}{2p_n}\right)\left(1-\left(t^\Omega_{*,n}(u_n)\right)^2\right)\|u_n\|^2_\lambda-\left(\frac{1}{2p_n}-\frac{1}{2\cdot 2^*_{\sigma,s}}\right)\left(t^\Omega_{*,n}(u_n)\right)^2\|u_n\|^2_\lambda.
	\end{align*}
	Similar to \eqref{<1} in the proof of Proposition \ref{=}, we have $t^\Omega_{*,n}(u_n)\le 1+o_n(1)$. Hence
	\begin{equation*}
	\left(\frac{1}{2}-\frac{1}{2p_n}\right)\left(1-\left(t^\Omega_{*,n}(u_n)\right)^2\right)\|u_n\|^2_\lambda\ge o_n(1),
	\end{equation*}
	and by $p_n\to2^*_{\sigma,s}$ as $n\to\infty$, we have
	\begin{equation*}
	\left(\frac{1}{2p_n}-\frac{1}{2\cdot 2^*_{\sigma,s}}\right)\left(t^\Omega_{*,n}(u_n)\right)^2\|u_n\|^2_\lambda=o_n(1).
	\end{equation*}
	Therefore,
	\begin{equation*}
	J_{\varepsilon_n}(u_n)-J^\Omega_*\left(t^\Omega_{*,n}(u_n)u_n\right)\ge o_n(1).
	\end{equation*}
	Due to \eqref{nto}, we get
	\begin{equation*}
	\vartheta^\Omega_*\le\lim_{n\to\infty}J^\Omega_*\left(t^\Omega_{*,n}(u_n)u_n\right)\le\lim_{n\to\infty}\left(J_{\varepsilon_n}(u_n)+o_n(1)\right)=\vartheta_*=\vartheta^\Omega_*,
	\end{equation*}
	which implies 
	\begin{equation*}
    \lim_{n\to\infty}J^\Omega_*\left(t^\Omega_{*,n}(u_n)u_n\right)=\vartheta^\Omega_*.
	\end{equation*}
	In terms of Ekeland variational principle (see Theorem 8.5 in Willem \cite{Willem}), there exist sequences $\{v_n\}_{n\in\mathbb N}\subset\mathcal N^\Omega_*$ and $\{\nu_n\}_{n\in\mathbb N}\subset\mathbb R$ such that as $n\to\infty$,
	\begin{equation}\label{approximate}
	\|v_n-t^\Omega_{*,n}(u_n)u_n\|_\lambda\to0,
	\end{equation}
	\begin{equation*}
	J^\Omega_*(v_n)=\left(\frac{1}{2}-\frac{1}{2\cdot 2^*_{\sigma,s}}\right)\|v_n\|^2_\lambda\to\vartheta^\Omega_*,
	\end{equation*}
	\begin{equation*}
	(J^\Omega_*)'(v_n)-\nu_n(G^\Omega_*)'(v_n)\to0\quad\text{in }  H^{-s}_0(\Omega),
	\end{equation*}
	where $G^\Omega_*(v_n):=\left<(J^\Omega_*)'(v_n),v_n\right>$. By Lemma \ref{PS}, we get that $\{v_n\}_{n\in\mathbb N}$ is a (PS)-sequence for the free functional $J^\Omega_*$ at level $\vartheta^\Omega_*$. Then Remark \ref{k=1} implies that 
	\begin{equation*}
	v_n-(S^s_{H,L})^{\frac{1}{2\cdot2^*_{\sigma,s}-2}}U_{R_n,x_n}\to0\quad\text{in } D^{s,2}(\mathbb R^N) \ \text{ as }n\to\infty,
	\end{equation*}
	where $\{x_n\}_{n\in\mathbb N}\subset\Omega$ and $R_n\to+\infty$ as $n\to\infty$. Write
	\begin{equation*}
	v_n=(S^s_{H,L})^{\frac{1}{2\cdot2^*_{\sigma,s}-2}}U_{R_n,x_n}+w_n, 
	\end{equation*}
    where $\|w_n\|_{s,2}\to0$ as $ n\to\infty$, which implies $|w_n|_{2^*_s}\to0$ as $ n\to\infty$. Unless to relabel $w_n$, we get from \eqref{approximate} that
	\begin{equation*}
	t^\Omega_{*,n}(u_n)u_n=(S^s_{H,L})^{\frac{1}{2\cdot2^*_{\sigma,s}-2}}U_{R_n,x_n}+w_n.
	\end{equation*}
    According to \eqref{barycenter}, we get that for $x=(x^1,x^2,\cdots,x^N)\in\mathbb R^N$,
	\begin{equation}\label{3terms}
	\begin{aligned}
	&\beta^i\left(t^\Omega_{*,n}(u_n)u_n\right)\big|t^\Omega_{*,n}(u_n)u_n\big|^{2^*_s}_{2^*_s}=\int_{\mathbb R^N}x^i\big|t^\Omega_{*,n}(u_n)u_n\big|^{2^*_s}\text{d}x\\
	=&\int_{\mathbb R^N}x^i\Big(\big|(S^s_{H,L})^{\frac{1}{2\cdot2^*_{\sigma,s}-2}}U_{R_n,x_n}+w_n\big|^{2^*_s}-\big|(S^s_{H,L})^{\frac{1}{2\cdot2^*_{\sigma,s}-2}}U_{R_n,x_n}\big|^{2^*_s}\Big)\text{d}x\\
	&+(S^s_{H,L})^{\frac{2^*_s}{2\cdot2^*_{\sigma,s}-2}}\int_{\mathbb R^N}x^i|U_{R_n,x_n}|^{2^*_s}\text{d}x\\
	=:&I_1+I_2.
	\end{aligned}
	\end{equation}
					 		
	Since $U_{R_n,x_n}(x)=R_n^{\frac{N-2s}{2}}U_{1,0}(R_n(x-x_n))$, we take the change of variables $\tilde x=R_n(x-x_n)$ and use the property of integral for odd functions in symmetric domain to obtain that
	\begin{equation}\label{I2}
	\begin{aligned}
	I_2&=(S^s_{H,L})^{\frac{2^*_s}{2\cdot2^*_{\sigma,s}-2}}\Big(x^i_n|U_{1,0}|^{2^*_s}_{2^*_s}+\frac{1}{R_n}\int_{\mathbb R^N}\tilde x^i|U_{1,0}(\tilde x)|^{2^*_s}\text{d}\tilde x\Big)\\
	&=(S^s_{H,L})^{\frac{2^*_s}{2\cdot2^*_{\sigma,s}-2}}x^i_n|U_{1,0}|^{2^*_s}_{2^*_s}.
	\end{aligned}
	\end{equation}
	
	Notice that $\{v_n\}_{n\in\mathbb R^N}$ is supported in $\Omega$, we have 
	\begin{align*}
	w_n=-(S^s_{H,L})^{\frac{1}{2\cdot2^*_{\sigma,s}-2}}U_{R_n,x_n}\quad\text{in }\Omega^c.
	\end{align*}			
	By the mean value theorem, the H\"older inequality and  the change of variables $\tilde x=R_n(x-x_n)$, we get that
	\begin{equation}\label{I1}
	\begin{aligned}
	|I_1|\le&2^*_s\int_{\mathbb R^N}|x^i|\big|(S^s_{H,L})^{\frac{1}{2\cdot2^*_{\sigma,s}-2}}U_{R_n,x_n}+\theta w_n\big|^{2^*_s-1}|w_n|\text{d}x \ (0<\theta<1)\\
	\le&C\int_\Omega\big(|U_{R_n,x_n}|^{2^*_s-1}|w_n|+|w_n|^{2^*_s}\big)\text{d}x+C\int_{\Omega^c}|x^i||U_{R_n,x_n}|^{2^*_s}\text{d}x\\
	\le&C\big(|U_{1,0}|^{2^*_s-1}_{2^*_s}|w_n|_{2^*_s}+|w_n|^{2^*_s}_{2^*_s}\big)+C\int_{R_n(\Omega^c-x_n)}\frac{\frac{1}{R_n}|\tilde x^i|+|x^i_n|}{(1+|\tilde x|^2)^N}\text{d}\tilde x\\
	=&o_n(1).
	\end{aligned}
    \end{equation}

	Similarly, we also have
	\begin{equation}\label{I0}
	\begin{aligned}
	\big|t^\Omega_{*,n}(u_n)u_n\big|^{2^*_s}_{2^*_s}=(S^s_{H,L})^{\frac{2^*_s}{2\cdot2^*_{\sigma,s}-2}}|U_{1,0}|^{2^*_s}_{2^*_s}+o_n(1).
	\end{aligned}
	\end{equation}
	Inserting \eqref{I2}, \eqref{I1} and \eqref{I0} into \eqref{3terms}, we obtain that
	\begin{align}\label{I123}
	\beta^i(u_n)=\beta^i\left(t^\Omega_{*,n}(u_n)u_n\right)=\frac{x^i_n(S^s_{H,L})^{\frac{N}{N-\sigma+2s}}\big|U_{1,0}\big|^{2^*_s}_{2^*_s}+o_n(1)}{(S^s_{H,L})^{\frac{N}{N-\sigma+2s}}\big|U_{1,0}\big|^{2^*_s}_{2^*_s}+o_n(1)}.
	\end{align}
	Noticing $\{x_n\}_{n\in\mathbb N}\subset\Omega$, we get from \eqref{I123} that $\beta(u_n)\in\Omega^+_r$ for $n$ sufficiently large, which contradicts with \eqref{contradiction}. This completes the proof.
\end{proof}

Now we are ready to prove the main theorem in this paper.
\vskip0.1in
\textbf{\large Proof of Theorem \ref{t}:}
\vskip0.1in
\textbf{Step 1.} We show the existence of $cat(\bar\Omega)$ low energy solutions for the problem \eqref{target}.

Let's fix $\delta_0>0$ and $\varepsilon_0(\delta_0)>0$ as in Proposition \ref{beta}. Then for all $\varepsilon<\varepsilon_0(\delta_0)$, there holds
\begin{equation}\label{arrow}
    u\in\mathcal N_\varepsilon\ \text{ and }\ J_\varepsilon(u)<\vartheta_\varepsilon+\delta_0\Longrightarrow \beta(u)\in\Omega^+_r.
\end{equation}
By Proposition \ref{=}, for the above $\delta_0$,  there exists $\bar\varepsilon(\delta_0)$ such that
\begin{equation}\label{vartheta}
    |\vartheta_\varepsilon-\vartheta_*|\le\frac{\delta_0}{2},\quad \forall \varepsilon<\bar\varepsilon(\delta_0).
\end{equation}
Due to the proof of Lemma \ref{<}, there exists $\tilde\varepsilon(\delta_0)>0$ such that for all $\varepsilon<\tilde\varepsilon(\delta_0)$, there is $R=R(\delta_0,\varepsilon)>1$ such that
\begin{equation}\label{J}
    J_\varepsilon(t_\varepsilon(u_{R,x_0})u_{R,x_0}(x))\le\vartheta_*+\frac{\delta_0}{2},
\end{equation}
where $u_{R,x_0}$ is defined in \eqref{u} and $t_\varepsilon(u_{R,x_0})>0$ is the unique value satisfying $t_\varepsilon(u_{R,x_0})u_{R,x_0}\in\mathcal N_\varepsilon$.

After taking $0<\varepsilon<\min\{\varepsilon_0(\delta_0),\bar\varepsilon_0(\delta_0),\tilde\varepsilon_0(\delta_0)\}$ and choosing $R=R(\delta_0,\varepsilon)>1$ sufficiently large, we define
\begin{align}\label{phi}
    \varphi_\varepsilon:\Omega^-_r\to\mathcal N_\varepsilon \ \text{ and }\ \varphi_\varepsilon(x):=t_\varepsilon(u_{R,x_0})u_{R,x_0}(x).
\end{align}
It follows from \eqref{vartheta}, \eqref{J} and \eqref{phi} that
\begin{align}\label{subset}
    \varphi_\varepsilon(\Omega^-_r)\subseteq\mathcal N_\varepsilon\cap J_\varepsilon^{\vartheta_*+\delta_0/2}\subseteq\mathcal N_\varepsilon\cap J_\varepsilon^{\vartheta_\varepsilon+\delta_0},
\end{align}
where $J_\varepsilon^c:=\{u\in H^s_0(\Omega):J_\varepsilon(u)\le c\}$ $(c\in\mathbb R)$ denotes the level set of $J_\varepsilon$.

By means of \eqref{arrow} and \eqref{subset}, the following maps are well-defined
\begin{equation*}
    \Omega^-_r\stackrel{\varphi_\varepsilon}{\longrightarrow}\mathcal N_\varepsilon\cap J_\varepsilon^{\vartheta_\varepsilon+\delta_0}\stackrel{\beta}{\longrightarrow}\Omega^+_r,
\end{equation*}
and $\beta\circ\varphi_\varepsilon$ is homotopic to the identity on $\Omega^-_r$. Due to Remark \ref{BC}, we get that
\begin{equation*}
    cat_{J_\varepsilon^{\vartheta_\varepsilon+\delta_0}}(\varphi_\varepsilon(\Omega^-_r))\ge cat(\Omega^-_r)=cat(\bar\Omega)>1.
\end{equation*}
Thus we find a sublevel of $J_\varepsilon$ on $\mathcal N_\varepsilon$ with category greater than $cat(\bar\Omega)$. Notice that $J_\varepsilon$ satisfies the (PS)-condition on $\mathcal N_\varepsilon$, we obtain from Proposition \ref{BCP} that there exist at least $cat(\bar\Omega)$ critical points of $J_\varepsilon$ in $J_\varepsilon^{\vartheta_\varepsilon+\delta_0}$, which correspond to the low energy solutions of the problem \eqref{target}.

\textbf{Step 2.}  We prove the existence of another high energy solution for \eqref{target} when $\Omega$ is not contractible.

Given arbitrarily a positive function $v\in D^{s,2}(\mathbb R^N)$ and $x_0\in\Omega^-_r$, we let
\begin{equation*}
    \bar v(x):=v(x)\phi_{r,x_0}(x),
\end{equation*}
where $\phi_{r,x_0}$ is defined as in \eqref{cutoff}. Then $\bar v\in H^s_0(\Omega)$. Set 
\begin{equation*}
    E_\varepsilon:=\left\{\theta\bar v(x)+(1-\theta)e(x):\theta\in[0,1],e\in \varphi_\varepsilon(\Omega^-_r)\right\},
\end{equation*}
then $\varphi_\varepsilon(\Omega^-_r)\subset E_\varepsilon\subset H^s_0(\Omega)$. Moreover, $E_\varepsilon$ is compact and contractible in $H^s_0(\Omega)$. Notice that the functions in $\varphi_\varepsilon(\Omega^-_r)$ are all positive, we get that $E_\varepsilon$ contains only positive functions. For $u\in E_\varepsilon$, there is a unique $t_\varepsilon(u)>0$ such that $t_\varepsilon(u)u\in\mathcal N_\varepsilon$, and then
\begin{equation}\label{tu=}
    (t_\varepsilon(u))^{2p_\varepsilon-2}=\frac{\|u\|^2_\lambda}{\displaystyle\int_{\Omega\times\Omega}\frac{(u^+(x))^{p_\varepsilon}(u^+(y))^{p_\varepsilon}}{|x-y|^\sigma}\text{d}x\text{d}y}=\frac{\|u\|^2_\lambda}{\displaystyle\int_{\Omega\times\Omega}\frac{u^{p_\varepsilon}(x)u^{p_\varepsilon}(y)}{|x-y|^\sigma}\text{d}x\text{d}y}.
\end{equation}
Define
\begin{equation*}
    T_\varepsilon:=\left\{t_\varepsilon(u)u:u\in E_\varepsilon\right\},
\end{equation*}
then $\varphi_\varepsilon(\Omega^-_r)\subset T_\varepsilon\subset\mathcal N_\varepsilon$. Additionally, $T_\varepsilon$ is compact and contractible in $\mathcal N_\varepsilon$, and $T_\varepsilon$ contains only positive functions. Finally, we denote
\begin{equation*}
    m_\varepsilon:=\max_{u\in E_\varepsilon} J_\varepsilon(t_\varepsilon(u)u)=\max_{u\in T_\varepsilon}J_\varepsilon(u).
\end{equation*}
Then $T_\varepsilon\subset\mathcal N_\varepsilon\cap J_\varepsilon^{m_\varepsilon}$ and $m_\varepsilon\ge\vartheta_\varepsilon$.

\textbf{Claim:} There exists a constant $c>0$ such that for each $\varepsilon>0$ small, it holds $m_\varepsilon<c$.

\begin{proof}[\bf Proof]
	For $u\in E_\varepsilon$, we have
	\begin{equation*}
	J_\varepsilon(t_\varepsilon(u)u)=\frac{p_\varepsilon-1}{2p_\varepsilon}(t_\varepsilon(u))^2\|u\|^2_\lambda.
	\end{equation*}
	
	It follows from the definition of $E_\varepsilon$ and \eqref{phi}, \eqref{tofu} that
	\begin{equation}\label{s2<C}
	\|u\|_{s,2}\le\|\bar v\|_{s,2}+2(S^s_{H,L})^{\frac{1}{2\cdot2^*_{\sigma,s}-2}}\|u_{R,x_0}\|_{s,2}\le C
	\end{equation}
	and 
	\begin{equation}\label{2<C}
	|u|_2\le|\bar v|_2+2(S^s_{H,L})^{\frac{1}{2\cdot2^*_{\sigma,s}-2}}|u_{R,x_0}|_2\le C.
	\end{equation}
	Since the domain $\Omega$ is bounded, we denote with $\text{diam}\Omega$ as its diameter. For any $x,y\in\Omega$, we have $|x-y|\le 2\text{diam}\Omega$. Hence for $u\in E_\varepsilon$,
	\begin{equation}\label{>C}
	\int_{\Omega\times\Omega}\frac{u^{p_\varepsilon}(x)u^{p_\varepsilon}(y)}{|x-y|^\sigma}\text{d}x\text{d}y
	\ge\frac{1}{(2\text{diam}\Omega)^\sigma}\int_{B_{r/2}(x_0)\times B_{r/2}(x_0)}u^{p_\varepsilon}(x)u^{p_\varepsilon}(y)\text{d}x\text{d}y>C>0.
	\end{equation}
	Indeed, we get from \eqref{tofu} that for $\varepsilon>0$ small enough and $x\in B_{r/2}(x_0)$, 
	\begin{align*}
	u(x)\ge&\theta v(x)+(1-\theta)1/2(S^s_{H,L})^{\frac{1}{2\cdot2^*_{\sigma,s}-2}}U_{R,x_0}(x)\\
	\ge&\max\Big\{\theta v(x),(1-\theta)1/2(S^s_{H,L})^{\frac{1}{2\cdot2^*_{\sigma,s}-2}}U_{R,x_0}(x)\Big\}\\
	\ge&\max\Big\{\theta \min\Big\{v,1/2(S^s_{H,L})^{\frac{1}{2\cdot2^*_{\sigma,s}-2}}U_{R,x_0}\Big\},(1-\theta)\min\Big\{v,1/2(S^s_{H,L})^{\frac{1}{2\cdot2^*_{\sigma,s}-2}}U_{R,x_0}\Big\}\Big\}\\
	\ge&1/2\min\Big\{v,1/2(S^s_{H,L})^{\frac{1}{2\cdot2^*_{\sigma,s}-2}}U_{R,x_0}\Big\},\quad \forall u\in E_\varepsilon,
	\end{align*}
	and then $u^{p_\varepsilon}(x)\ge C\min\left\{v,U_{R,x_0}\right\}^{2^*_{\sigma,s}}$, which ensures the correctness of \eqref{>C}. 
	
	According to \eqref{tu=}, \eqref{s2<C}, \eqref{2<C} and \eqref{>C}, we see that $t_\varepsilon(u)$ is bounded on $E_\varepsilon$ uniformly in $\varepsilon$. Concequently, $J_\varepsilon$ is bounded on $T_\varepsilon$ uniformly in $\varepsilon$. This finishes the proof of the claim.	
\end{proof}
    
Similarly to Section 6 in Benci et al \cite{Benci-Bonanno-Micheletti} and as the same argument of Proposition \ref{BCP} in the contractible case, we conclude that there exists another solution $\hat u$ to the problem \eqref{target} such that $$\vartheta_\varepsilon+\delta_0<J_\varepsilon(\hat u)\le m_\varepsilon.$$
We complete the proof of Theorem \ref{t}.

\begin{center}
	{\bf{Acknowledgement}} 
\end{center}
    We would like to show our sincere gratitude to the reviewers for their valuable suggestions, and especially for the relevant references suggested.





\vskip0.7in

\footnotesize


\end{document}